\documentclass[11pt]{article}
\usepackage{amssymb, hyperref}
\newtheorem{thm}{Theorem}
\newtheorem{prop}{Proposition}
\newtheorem{cor}{Corollary}

\newtheorem{lemma}{Lemma}
\newtheorem{defn}{Definition}

\begin{document}

\parindent 0mm
\parskip 2.5mm

\setlength{\textheight}{22cm} 
\headsep=15pt
\setlength{\textwidth}{15.5cm}
\setlength{\oddsidemargin}{0.5cm} \setlength{\topmargin}{-.5cm}
\setlength{\evensidemargin}{\oddsidemargin}

\voffset=0.75cm
\hoffset=0.3cm
\vsize=23cm 
\hsize=15.5cm 
\mathsurround=2pt
\baselineskip 12.8pt
\hfuzz=3.5pt



\def\noheadline{\hfill}
\def\runningheadline{\ifodd\pageno\rightheadline \else\leftheadline\fi}
\def\rightheadline{\rmh\hfil{\slh
                Order $p$ automorphisms of the open disc of a
                $p\!$-adic field}\hfil\folio}
\def\leftheadline{\rmh\folio\hfil{\slh
   Order $p$ automorphisms of the open disc of a $p\!$-adic field}\hfil}
\baselineskip 12.8pt
\font\slh=cmsl10
\font\rmh=cmr10
\font\preloaded=eufm10\font\preloaded=eufm7\font\preloaded=eufm5
\font\preloaded=msbm10\font\preloaded=msbm7\font\preloaded=msbm5
\font\teneu=eufm10\font\seveneu=eufm7\font\fiveeu=eufm5
\font\tenlv=msbm10\font\sevenlv=msbm7\font\fivelv=msbm5



\def\eu #1{{\mathchoice{{\hbox{\teneu #1}}}{{\hbox{\teneu #1}}}
{{\hbox{\seveneu #1}}}{{\hbox{\fiveeu #1}}}}}
\def\lv #1{{\mathchoice{{\hbox{\tenlv #1}}}{{\hbox{\tenlv #1}}}
{{\hbox{\sevenlv #1}}}{{\hbox{\fivelv #1}}}}}


\def\cA{{\cal A}}\def\cB{{\cal B}}\def\cC{{\cal C}}\def\cD{{\cal D}}
\def\cE{{\cal E}}\def\cF{{\cal F}}\def\cG{{\cal G}}\def\cH{{\cal H}}
\def\cI{{\cal I}}\def\cJ{{\cal J}}\def\cK{{\cal K}}\def\cL{{\cal L}}
\def\cM{{\cal M}}\def\cN{{\cal N}}\def\cO{{\cal O}}\def\cP{{\cal P}}
\def\cQ{{\cal Q}}\def\cR{{\cal R}}\def\cS{{\cal S}}\def\cT{{\cal T}}
\def\cU{{\cal U}}\def\cV{{\cal V}}\def\cW{{\cal W}}\def\cX{{\cal X}}
\def\cY{{\cal Y}} \def\cZ{{\cal Z}}


\def\dbO{{\eu O}}\def\dbm{{\eu m}}\def\dbP{{\eu P}}\def\dbQ{{\eu Q}}
\def\dbp{{\eu p}}\def\dbq{{\eu q}}
\def\dbb{{\eu b}}\def\dbm{{\eu m}}
\def\dba{{\eu a}}


\def\lvA{{\lv A}}\def\lvC{{\lv C}}\def\lvF{{\lv F}}\def\lvG{{\lv G}}
\def\lvK{{\lv K}}\def\lvN{{\lv N}}\def\lvP{{\lv P}}\def\lvQ{{\lv Q}}
\def\lvR{{\lv R}}\def\lvZ{{\lv Z}}


\def\bftot#1{{\setbox0=\hbox{#1}\setbox1=\hbox{\hbox to.1pt{}\copy0}\copy1\kern-\wd0%
\copy1\kern-\wd0\copy1\kern-\wd0\copy1\kern-\wd0\copy1\kern-\wd0\copy1}}

\def\bfmatten#1{{\mathsurround=0pt\bftot{${\textstyle{#1}}$}}}
\def\bfmatseven#1{{\mathsurround=0pt\bftot{${\scriptstyle{#1}}$}}}
\def\bfmatfive#1{{\mathsurround=0pt\bftot{${\scriptscriptstyle{#1}}$}}}

\def\bfmath#1{{\mathchoice{{\hbox{\bfmatten{#1}}}}{{\hbox{\bfmatten{#1}}}}
{{\hbox{\bfmatseven{#1}}}}{{\hbox{\bfmatfive{#1}}}}}}

\font\secsize=cmbx10 scaled\magstep1
\font\capsize=cmbx10 scaled\magstep2 
\font\titlesize=cmbx10 scaled\magstep3

\def\bs{\bigskip}
\def\ms{\medskip}
\def\sms{\smallskip}
\def\bsn{\bigskip\noindent}
\def\msn{\medskip\noindent}
\def\ssn{\smallskip\noindent}
\def\pn{\par \noindent}
\def\hb{\hfill\break}

\def\commdiag#1#2#3#4#5#6#7#8{
              $$\def\normalbaselines{
                    \baselineskip20pt\lineskip3pt\lineskiplimit3pt}
              \matrix{#1&\lghor{#5}&#2&\cr
                      \dwn{#6}&&\dwn{#7}&\cr
                      #3&\lghor{#8}&#4&\cr}
                              $$}
    \def\congg{\hhb1\cong\hhb1}
\def\commdiagup#1#2#3#4#5#6#7#8{
              $$\def\normalbaselines{
                    \baselineskip20pt\lineskip3pt\lineskiplimit3pt}
              \matrix{#1&\lghor{#5}&#2&\cr
                      \uparrow{#6}&&\uparrow{#7}&\cr
                      #3&\lghor{#8}&#4&\cr}
                              $$}
    \def\congg{\hhb1\cong\hhb1}

\def\updiag#1#2#3#4{
              $$\def\normalbaselines{
                    \baselineskip20pt\lineskip3pt\lineskiplimit3pt}
              \matrix{#1&\lgrghtar&#2&\cr
                      \Big\uparrow&&\Big\uparrow&\cr
                      #3&\lgrghtar&#4&\cr}
                              $$}

\def\diag#1#2#3#4{
             $$\def\normalbaselines{
                          \baselineskip20pt\lineskip3pt\lineskiplimit3pt}
                \matrix{#1&\lgrghtar&#2&\cr
                         \Big\downarrow&&\Big\downarrow&\cr
                        #3&\lgrghtar&#4&\cr}
                                 $$}

\def\qtdiag#1#2#3{
                   $$\def\normalbaselines{
                          \baselineskip20pt\lineskip3pt\lineskiplimit3pt}
                          \matrix{#1&\lgrghtar&#2&\cr
                                & \Big\downarrow & \Big\downarrow&\cr
                                &&#3&\cr}
                                   $$}
\def\GQ{\Gal(\overline{\lvQ} /\lvQ)} 


\def\diagdiag#1#2#3#4#5#6#7#8{
              $$\def\normalbaselines{
                    \baselineskip20pt\lineskip3pt\lineskiplimit3pt}
              \matrix{#1&#5&#2&\cr
                      #6&  &#7&\cr
                      #3&#8&#4&\cr}
                              $$}

\def\Chabdiag#1#2#3#4#5#6#7#8#9{
              $$\def\normalbaselines{
                    \baselineskip20pt\lineskip3pt\lineskiplimit3pt}
              \matrix{#1&#5&#2&&\cr
                      #6&  &#7&&\cr
                      #3&#8&#4&#9&\cr}
                             $$}
\def\diagram#1{$$\def\normalbaselines{
                    \baselineskip20pt\lineskip3pt\lineskiplimit3pt}
                    #1
                    $$}

 \def\drpr#1#2{{\partial{#1}\over\partial{#2}}}
 \def\drv#1#2#3{{\displaystyle\biggl|\drpr{#1}{#2}{#3}\biggr|}}
 \def\dwn#1{\Big\downarrow\rlap{$\vcenter{\hbox{$\scriptstyle{#1}$}}$}}
 \def\upa#1{\Big\uparrow\rlap{$\vcenter{\hbox{$\scriptstyle{#1}$}}$}}
 \def\hhb#1{\hbox to#1pt{}}
 \def\hor#1{\smash
          {\mathop{{\lgrghtar}}\limits^{\lower2pt\hbox{$\scriptstyle{#1}$}}}}
 \def\lghor#1{\smash
         {\mathop{{\lglgrghtar}}\limits^{\lower2pt\hbox{$\scriptstyle{#1}$}}}}
 \def\llhor#1{\smash
         {\mathop{{\lgleftar}}\limits^{\lower2pt\hbox{$\scriptstyle{#1}$}}}}

 \def\lhor#1{\smash
         {\mathop{{\lglgleftar}}\limits^{\lower2pt\hbox{$\scriptstyle{#1}$}}}}
\def\binj#1{\smash
         {\mathop{{\inj}}\limits^{\lower2pt\hbox{$\scriptstyle{#1}$}}}}

\def\inj{\hookrightarrow}
\def\linj{\hookleftarrow}
\def\lglgrghtar{{\relbar\joinrel\relbar\joinrel\rightarrow}}
\def\lgrghtar{{\relbar\joinrel\rightarrow}}
\def\lgleftar{{\leftarrow\joinrel\relbar}}
\def\lglgleftar{{\leftarrow\joinrel\relbar\joinrel\relbar}}
\def\lra{{\longrightarrow}}
\def\lla{{\longleftarrow}}
\def\upa#1{\Big\uparrow\rlap{$\vcenter{\hbox{$\scriptstyle{#1}$}}$}}
\def\rdratmap{\hbox{\kern-6pt\hbox{$-\mkern1.5mu-\mkern-2mu$}
            \raise0.8pt\hbox{$\scriptstyle\rightarrow$}\kern-6pt}}
\def\rratmap{\hbox{\kern-0pt\hbox{$-\mkern1.5mu-$}
            \kern-5pt\raise0.8pt\hbox{$\scriptstyle\rightarrow$}\kern-0pt}}
\def\dratmap{\vbox{\baselineskip=6pt \lineskiplimit=0pt
              \kern5pt{\hbox{$\mapstochar$}\hbox{$\mapstochar$}
              \kern1pt\hbox{$\mkern-2.3mu\scriptstyle\downarrow$}}}}
\def\sag#1{\hhb1{\hbox to#1mm{\rightarrowfill}}\hhb1}
\def\sdp{\hbox to10pt{\hss\hbox{\mathsurround=0pt$\times$\kern-1.6pt
       \hbox{\vrule height4.5pt width.6pt}\hbox to1.6pt{}}\hss}}

 \def\formel#1{ $$\abovedisplayskip=.4\abovedisplayskip
                  \belowdisplayskip=.4\belowdisplayskip
                                   #1 $$}
\def\qed{\hhb2\vbox{\hrule\hbox{\vrule\kern 2pt\vbox{\kern 6pt}
                                               \kern2pt\vrule}\hrule}}
\def\vid{{\hhb{-3.25}\not\hhb{-2.25}\lower-1.25pt
     \hbox{\mathsurround=0pt$\scriptscriptstyle\bigcirc$}}}
\def\bcul#1{{\scriptstyle\bigcup\limits_{{#1}}\hhb{1.3}}}
\def\bcull#1#2{{\scriptstyle\bigcup\limits_{{#1}}^{{#2}}\hhb{1.3}}}
\def\bcal#1{{\scriptstyle\bigcap\limits_{{#1}}\hhb{1.3}}}
\def\plim{\lim\limits_{\lower-2pt\hbox{$\longleftarrow$}}}
\def\prodl#1{{\textstyle\prod\limits_{#1}\hhb1}}
\def\prodll#1#2{{\textstyle\prod\limits_{#1}^{#2}\hhb1}}
\def\suml#1{{\textstyle\sum\limits_{#1}\hhb1}}
\def\sumll#1#2{{\textstyle\sum\limits_{#1}^{#2}\hhb1}}
\def\proof{{\it Proof}}
\def\iso{\approx}
\def\sfrac#1#2{{ \textstyle{{{#1}}\over{{#2}}} }}  
\def\deg{{\rm deg}\,}
\def\dim{{\rm dim}}
\def\mod{{\rm mod}\, }
\def\bua{\big\uparrow}
\def\lra{\longrightarrow}
\def\Eul#1#2{{\chi(#1\vert#2)}}
\def\frac#1#2{{\displaystyle{{#1}\over{#2}}}}
\def\sst{\scriptscriptstyle\times}
\def\id{{{\bf id}}}
\def\bmmu{{\bfmath{\mu}}}
\def\G{{\lvG}}
\def\nltimes{{\hbox{$\raise 1.3pt\hbox{$\scriptscriptstyle |\!$}
              \!\!\times$}}}
\def\ro{{\rho}}
\def\el{{\ell}}
\def\Spec{{\rm Spec}\, }
\def\Gal{{\rm Gal}\, }
\def\Cov{{\rm Cov}\, }
\def\End{{\rm End}\, }
\def\Aut{{\rm Aut}\, }
\def\Im{{\rm Im}\, }
\def\Re{{\rm Re}\, }
\def\Res{{\rm Res}\, }
\def\La{{\Lambda}}
\def\gep{{\gamma_{1, \varepsilon}}}
\def\Log{{\rm Log}\, }
\def\Li{{\rm Li}\, }
\def\Rlim{{\rm Rlim}\, }
\def\arg{{\rm arg}\, }
\def\exp{{\rm exp}\, }
\def\Isom{{\rm Isom}\, }
\def\ord{{\rm ord}\, }
\def\Pomt{{\lvP^{1} \backslash \{0,1,\infty\}}}

{\title{On a generalization of {\sc Chen's} iterated integrals}}
{\author{Sheldon Joyner}}

\maketitle

{\sc The University of Western Ontario
\newline
Department of Mathematics
\newline 
Middlesex College
\newline London
\newline Ontario N6A 5B7
\newline Canada}
\newline
{\em  email:} sjoyner at uwo dot ca
\newline
{\em Fax:} 1-519-661-3610

\thispagestyle{empty}

\newpage
{\abstract{In this paper, 
{\sc Chen's} iterated integrals are generalized by interpolation of functions of
the positive integer number of times which particular forms are iterated in
integrals along specific paths, to certain complex values. These generalized
iterated integrals satisfy both an additive 
iterative property and comultiplication formula. In a particular example, a (non-classical) multiplicative iterative property is also shown to hold.
After developing this theory in the first part of the paper we discuss various applications, including the expression of certain zeta functions as complex iterated integrals (from which an obstruction to the existence of a contour integration proof of the functional equation for the {\sc Dedekind} zeta function emerges); a way of thinking about complex iterated derivatives arising from a reformulation of a result of {\sc Gel'fand} and {\sc Shilov} in the theory of distributions; and a direct topological proof of the monodromy of polylogarithms. 
{}
}}

\newpage

{\section{Introduction}}
The iterated integrals of {\sc K.-T. Chen} arise in arithmetic situations, a famous
example of which is the occurrence of the polyzeta values (also called multiple zeta values)
as periods relating two distinct rational structures on the mixed {\sc Hodge} structure 
which comprises the {\sc Hodge} realization of 
the motivic fundamental group of $\lvP^{1} \backslash \{0,1,\infty\}$ with 
tangential base-point $\overrightarrow{01}$. In this paper, it is shown that more general objects, 
including the polyzeta functions themselves, may be viewed as iterated integrals 
of a sort generalizing the notion introduced by {\sc Chen},
and the eventual hope is that such objects could thereby also acquire further 
arithmetic significance.

A very general formulation of these iterated integrals is presented in the first section of the paper, in which it is shown that formal generalizations of the antipode and product formulas satisfied by {\sc Chen's} integrals may be ideated and then exploited to define complex iterated integrals. In particular, whenever the relevant integrals converge, then for differential 1-forms $\alpha$ and $\beta$ on some differential manifold $M$ on which $\gamma$ is a piece-wise smooth path, we define
\[
\int_{\gamma}\alpha \beta^{s-1}:=
\int_{0}^{1}\frac{1}{\Gamma(s)}\left(\int_{z}^{1}\gamma^{*}\beta\right)^{s-1}\gamma^{*}\alpha (z)
\]
where $z$ is a parameter on $[0,1]$ for the pullback of $\alpha$ under $\gamma,$ and $s$ is some complex number. Using the classical fact that the beta function has an expression as a quotient of values of the gamma function, the necessary iterative property is then established. 

This definition admits of the proof of a comultiplication formula, which is the first of the main results presented in the paper. This formula, which  is subject to certain technical conditions ensuring convergence of the sum, is given by
\[
\int_{\gamma \delta}\alpha\beta^{s} 
=
\int_{\delta}\alpha \beta^{s}
+
\sum_{n=0}^{\infty}\int_{\gamma}\alpha \beta^{n}\int_{\delta}\beta^{s-n}
\]
where $\gamma$ and $\delta$ are paths which may be concatenated, and 
\begin{equation}
\label{Weird0}
\int_{\delta}\beta^{s-n}
\end{equation}
is interpreted as 
\[
\left( \matrix{s \cr n} \right)
\cdot \frac{\int_{\delta}\beta^{s}}{\int_{\delta}\beta^{n}}
\]
for those $n$ for which  (\ref{Weird0}) does not converge. (See Theorem \ref{cOmUlt} of section \ref{CoMult}.)

This formula may be used to show that under certain conditions, the complex iterated integrals are homotopy functionals.

A further application of the formula is another of the principal results of the paper, namely a direct proof of the monodromy of the polylogarithm functions. This is preferable to the classical proof entailing use of {\sc Jonqui\`ere's} formula. These ideas are discussed in the last section of the paper, \ref{Mop}.

Before getting to this, we develop the theory in a most interesting example, namely that of $M = \Pomt.$ In this case, taking $\beta = \frac{dz}{z},$ the coincidence of our notion along paths $\gamma = [t,1]$ with the classical fractional integral is shown, thereby providing intrinsic geometric motivation for the fractional integral (along with the {\sc Mellin} transform).  On the other hand, it emerges that the necessary iterative property was known classically in this particular example. We demonstrate that this iterative property characterizes the complex iterated integrals in the case of $\beta = \frac{dz}{1-z}$ integrated over $[0,1].$ 


Beyond this, it is also possible to extend the formalism to multiple versions of complex iterated integrals, and we carry this out in order to subsequently express the polyzeta {\em functions} as iterated integrals.

Also, a non-classical multiplicative iterative property arises. This has an amusing consequence for the {\sc Riemann} zeta function, which is shown to admit a complex iterated integral expression
\[
\zeta(s)
=
\int_{[0,1]}\frac{dz}{1-z}\left(\frac{dz}{z}\right)^{s-1}
\]
that corresponds to {\sc Abel's} integral
\[
\zeta(s) = \frac{1}{\Gamma(s)}\int_{0}^{\infty}\frac{x^{s}}{e^{x}-1}\frac{dx}{x} \;.
\]
Using the multiplicative iterativity, we arrive at a family of expressions for $\zeta(s)$ indexed by positive integers $k$, in which the integral corresponding to $k=2$ is nothing other than the theta function integral which is the basis of the {\sc Fourier} analysis proof of the functional equation for $\zeta(s).$ 

It is also noteworthy that $\zeta(s)$ may be regarded as an integral transform of the rational function $\frac{z}{1-z},$ in keeping with the general philosophy that zeta functions should be rational. A similar statement holds for the {\sc Dirichlet} $L$-functions. However, an interesting difference is that while $\frac{z}{1-z}$ has a pole at $z=1,$ the rational function corresponding to $L(s, \chi)$ for non-trivial {\sc Dirichlet} character $\chi$ is non-singular at $z=1.$ On the other hand, $\zeta(s)$ is singular at $s=1,$ while $L(s, \chi)$ has no pole there. Although there is {\em a priori} no connection between the $z$ and $s$ coordinates, this correspondence turns out to hold quite generally. In fact, we prove in Theorem \ref{sz1} of section \ref{sec4} that if $F(z)$ satisfies a suitable boundedness condition (guaranteeing convergence of the relevant integral), and is meromorphic in a neighborhood of $z=1,$ then 
\[
L(F)(s)
:=
\int_{[0,1]}F(z)\frac{dz}{z}
\left(\frac{dz}{z}\right)^{s-1} 
\]
has a pole at $s=1$ if and only if $\frac{F(z)}{z}$ has a pole of non-zero residue at $z=1;$ and should $L(F)(s)$ have a pole at $s=1$, the residue is a sum of coefficients of the {\sc Laurent} series expansion of $F$ about $z=1.$ A consequence of this theorem is that the function $F_{K}(z)$ associated by means of a complex iterated integral expression to the {\sc Dedekind} zeta function of any number field for which the residue $\rho_{K}$ of the pole at $s=1$ is irrational (which is expected to hold for all number fields other than $\lvQ$), is not meromorphic near $z=1$ and hence not rational. Consequently, by a result due to {\sc Fatou} (see [9]), we can say that the function $F_{K}(z)$ is not even algebraic. Moreover, by a theorem of {\sc Petersson} (quoted in [2]), it follows that this function is not analytically continuable beyond the unit disc. 

In this way, we arrive at Corollary \ref{MainCor} of section \ref{sec4}, which is one of our principal results: Irrationality of $\rho_{K}$ is an obstruction to the existence of a proof of the analytic continuation and functional equation for the {\sc Dedekind} zeta function $\zeta_{K}(s)$ using the contour integral approach of {\sc Riemann's} first proof of the corresponding facts for $\zeta(s).$

The formalism  has proven useful in explaining a well-known result of {\sc Gel'fand} and {\sc Shilov} (presented in [10]I.\S 3.5) to the effect that the generalized function 
\[
\frac{x_{+}^{s-1}}{\Gamma(s)}
\]
admits an analytic continuation which at the negative integer $-n$ is the same as the $n$th derivative {\sc Dirac} measure $\delta^{(n)}$ - i.e. 
the value of a test function $\phi(x)$ against the generalized function
\[
\frac{x_{+}^{s-1}}{\Gamma(s)}|_{s=-n}
\]
over the  reals is given by 
\[
\left( -\frac{dx}{x} \right)^{n}\phi(x)|_{x=0}.
\]
In terms of complex iterated integrals (via a change of variables) this can be reformulated as in Theorem \ref{thm3} of section \ref{sec2}, which in essence is the statement that for suitable $F(z)$, then with notation as before, $L(F)(s)$ admits an analytic continuation with poles at most at a finite set of positive integers, which has values at
negative integers $-k$ given by
\[
\left( z\frac{d}{dz} \right) ^{k}F(z) |_{z=1}\;.
\]
This new perspective shows that we should think about
the differential operator
\[
\left(z \frac{d}{dz}\right)^{t}(\cdot) |_{z=1}
\]
as the analytic continuation of 
\[
\int_{[0,1]}(\cdot) \left(\frac{dz}{z}\right)^{s}
\]
to $s=-t.$

{\sc Riemann's} integral expression for the analytic continuation of $\zeta(s)$ may be modified to give a proof. 

Altering this proof in turn, the remarkable fact emerges that for any $w \in (0,1),$
\[
\int_{[w,1]}F(z)\frac{dz}{z}\left(\frac{dz}{z}\right)^{s-1}
\]
has the same analytic continuation to negative integers as does $L(F)(s)$.

The author would like to thank 
Professor {\sc Minhyong Kim} for 
his continued encouragement and 
valuable suggestions. He would also like to thank the referee and Professor {\sc J. Lagarias}, both of whom alerted him to inaccuracies in earlier versions of the paper. The referee also made several much-appreciated suggestions which have vastly clarified the exposition.

{\section{Iterated integrals along paths on complex manifolds}}
\label{gene}

Suppose throughout that $\alpha$ and $\beta$ are holomorphic 1-forms on a complex manifold 
$M$ and $\gamma$ is a piecewise smooth path in $M.$ 
The task at hand is to define 
\[
\int_{\gamma}\alpha \beta^{s-1}
\]
as an iterated integral, for suitable {\em complex} $s$.

This generalization should adhere to some kind of shuffle product 
generalizing the product on the {\sc Hopf} algebra of {\sc Chen's} iterated integrals.
The suitable form of this generalized product is not obvious, but 
repeated application of the usual shuffle product formula 
and use of a simple induction argument shows that  for any positive integer $n$, 
\begin{equation}
\label{it}
\left(\int_{\gamma}\beta \right)^{n} = n!\int_{\gamma}\beta^{n}.
\end{equation}
Here, the $n$-fold integration on the right side is reduced to 
a single integration on the left side. For example,
when $\gamma = [0,1]$ in $M=\lvC,$ geometrically this equation gives a transition 
between integration over the $n$-cube $[0,1]^{n}$ 
(the integral on the left side is an $n$-fold product of equal integrals, 
which by {\sc Fubini's} Theorem may be considered as a single integral over the cube), and 
integration over the time-ordered $n$-simplex 
\[
\{
(t_{1}, \ldots, t_{n})\in \lvR^{n}
|0\leq t_{1} \leq \ldots \leq  t_{n} \leq 1
\}
\]
(the integral on the right side, by the definition of iterated integrals). 
There are $n!$ such simplices which together form the $n$-cube and the permutation of the 
$t_{j}$ which shows this gives a change of variables yielding 
$n!$ equal integrals, the sum of which is the integral over the cube. 

The $n$ in (\ref{it}) can be interpolated to other complex arguments in an essentially unique way:
The gamma function is the unique function interpolating $n!$ 
having  certain nice properties (namely it satisfies the functional 
equation $\Gamma(x+1) = x\Gamma(x)$, has $\Gamma(1) = 1$, and when restricted to the positive 
reals has  convex logarithm). Moreover,
raising to the $n$th power is uniquely interpolated to other complex values via
\[
x^{s} = \exp(s\Log x),
\]
once a choice of the logarithm has been made, say 
$\Log z = \log|z|+i \arg (z) + 2 \pi i r$ for some $r \in \lvZ$ with $-\pi < \arg (z)<\pi$ 
(i.e. branch cut along the negative reals, or what is the same, $\Log$ has domain $\lvC \backslash \lvR_{<0}$).

The reason that this fact is significant is that in defining some 
kind of  complex power of the iterated integral - i.e. ascribing  meaning to 
integration against some object which gives a valid interpretation of 
complex power of a differential form - we have to somehow bypass integrating ``$s$ number of times'' 
for complex variable $s.$ 

Hence we can make the
{\defn{\label{Cii0}
For all $s$ for which the integral on the right hand side converges, we define
\[
\int_{\gamma}\beta^{s-1}\alpha
:=
\int_{0}^{1}\frac{1}{\Gamma(s)}\left(\int_{0}^{z}\gamma^{*}\beta \right)^{s-1}\gamma^{*}\alpha(z)
\]
where $z$ is a parameter on $[0,1]$ for $\gamma^{*}\alpha.$
}}

Recall now that the antipode property of iterated integrals is
\begin{equation}
\label{antip}
\int_{\gamma}\omega_{0} \ldots \omega_{r}
=
(-1)^{r+1}\int_{\gamma^{-1}}\omega_{r} \ldots \omega_{0}
\end{equation}
where $\gamma^{-1}$ is the inverse path to $\gamma$ defined by $\gamma^{-1}(t) = \gamma(1-t).$ 
Along with the shuffle product, we exploit the obvious analogue of this antipode property to define the type of integrals which will appear in the applications:
{\defn{
\label{Cii}
For all $s$ for which the integral on the right hand side converges, we define
\[
\int_{\gamma}\alpha \beta^{s-1}
:=
(-1)^{s}\int_{\gamma^{-1}}\beta^{s-1}\alpha\;.
\]}}
\newline Here the $(-1)^{s}$ factor is interpreted in a formal way by inverting the direction of integration, and the integral on the right side is as in Definition \ref{Cii0} - i.e. it is given by
\[
(-1)^{s}\int_{1}^{0}\frac{1}{\Gamma(s)}\left(\int_{1}^{z}\gamma^{*}\beta \right)^{s-1}\gamma^{*}\alpha(z)
:=
\int_{0}^{1}\frac{1}{\Gamma(s)}\left(\int_{z}^{1}\gamma^{*}\beta \right)^{s-1}\gamma^{*}\alpha(z)
\]
where $z$ is a parameter on $[0,1]$ for $\gamma^{*}\alpha.$

Synthesizing these definitions, we arrive at
{\defn{
\label{DualCii}
For those $(r,s) \in \lvC^{2}$ for which the integral on the right side converges, we define
\[
\int_{\gamma} \alpha^{r-1}\beta^{s}:=
\int_{[0,1]}\frac{1}{\Gamma(r)} \left[\left( \int_{0}^{z} \gamma^{*}\alpha\right)^{r-1} \gamma^{*}\beta \right] \cdot (\gamma^{*}\beta)^{s-1}.
\]
}}
\newline The resulting integral may now be interpreted in the light of Definition \ref{Cii} since $\left[\left( \int_{0}^{z} \gamma^{*}\alpha\right)^{r-1} \gamma^{*}\beta \right]$ is a 1-form. Doing so, we find
\begin{equation}
\label{dualcii}
\int_{\gamma} \alpha^{r-1}\beta^{s} = 
\frac{(-1)^{s}}{\Gamma(r)\Gamma(s)}\int_{1}^{0} \left( \int_{1}^{z}\gamma^{*}\beta\right)^{s-1} \left( \int_{0}^{z} \gamma^{*}\alpha \right)^{r-1} \gamma^{*}\beta(z).
\end{equation}
Applying the antipode property directly to $\int_{\gamma}\alpha^{r-1}\beta^{s}$ and then using Definition \ref{DualCii} itself, gives
\[
(-1)^{r+s-1}\int_{[1,0]}\frac{1}{\Gamma(s)}\left( \int_{1}^{z}\gamma^{*}\beta \right)^{s-1} \gamma^{*}\beta \cdot (\gamma^{*} \alpha)^{r-1}
=
\frac{(-1)^{r+s-1+r}}{\Gamma(s)\Gamma(r)}\int_{0}^{1}\left(\int_{0}^{z} \gamma^{*}\alpha \right)^{r-1} \left( \int_{1}^{z}\gamma^{*}\beta \right)^{s-1} \gamma^{*}\beta 
\]
which equals $\int_{\gamma}\alpha^{r-1} \beta^{s}$ via (\ref{dualcii}).

Any number of variations exist of the kind of integral here defined. For example, suppose that $\omega$ is a third holomorphic 1-form on $M$, and set
\begin{equation}
\label{newdef}
\int_{\gamma} \omega \beta^{s} \alpha := \int_{0}^{1} \left( \int_{[0,z]} 
\gamma^{*}\omega (\gamma^{*} \beta)^{s} \right) \gamma^{*} \alpha (z)
\end{equation}
for those $s$ for which the integral converges, 
where the integral over $[0,z]$ follows the pattern of Definition \ref{Cii}. Also, for suitable complex $r$ and $s$, take
\begin{equation}
\label{newdefn2}
\int_{\gamma}\alpha \beta^{r} \omega^{s-1}
:=
\frac{(-1)^{r+s}}{\Gamma(s)}
\int_{1}^{0}\left(\int_{1}^{t} \gamma^{*}\omega \right)^{s-1}\gamma^{*}\beta(\gamma^{*}\beta)^{r-1}\gamma^{*}\alpha\;.
\end{equation}
Here, (\ref{newdef}) should be used to interpret the right hand side, as $\left(\int_{1}^{t} \gamma^{*}\omega \right)^{s-1}\gamma^{*}\beta$ is a 1-form.

A well-definedness issue arises in Definition \ref{DualCii} in the case that $\beta$ and $\alpha$ are equal, and in (\ref{newdefn2}) should $\omega = \beta$. Clearing up the former problem is instrumental in showing that the latter is a non-issue, and thereby that the integrals of Definition \ref{Cii} are indeed iterated:
{\prop{{\bf [Iterative Property]}
\newline For suitable pairs $(v,w) \in \lvC^{2},$ in fact
\begin{equation}
\label{It}
\int_{\gamma}\alpha \beta^{v+w-1}= \int_{\gamma} (\alpha \beta^{v}) \beta^{w-1}.
\end{equation}
where the integral on the right side should be understood using (\ref{newdefn2}). Both sides of this equality are well-defined since for suitable $(v,u) \in \lvC,$ when $\alpha = \beta$, 
\begin{equation}
\label{SuffIt}
\int_{[t,1]}(\gamma^{*}\beta)^{v+u-1}:=
\frac{1}{\Gamma(v+u)}
\left(\int_{t}^{1}\gamma^{*}\beta\right)^{v+u-1}
=
\frac{1}{\Gamma(v)\Gamma(u)}\int_{t}^{1}
\left(\int_{t}^{z}\gamma^{*}\beta\right)^{v-1}\left(\int_{z}^{1}\gamma^{*}\beta\right)^{u-1}\gamma^{*}\beta
\end{equation}
is valid, where 
$z$ and 
$t$ are parameters on the interval - in particular we can take $u=1$ and $t=0.$}}
\newline {\bf Proof:} We begin by proving the well-definedness statement. To this end, define
\[
a:= \frac{\int_{t}^{z}\gamma^{*}\beta}{\int_{t}^{1}\gamma^{*}\beta} \; .
\]
Then notice that if $z=t,\; a=0$ whereas $a=1$ when $z=1.$ Also, 
\[
1-a = \frac{\int_{z}^{1}\gamma^{*}\beta}{\int_{t}^{1}\gamma^{*}\beta} \; .
\]
Finally, viewing $a$ as a function of $z$, 
\[
da = \frac{\gamma^{*}\beta(z)}{\left(\int_{t}^{1}\gamma^{*}\beta\right)} \; ,
\]
where $\gamma^{*}\beta(z)$ indicates that $z$ is the parameter of integration along $[0,1]$ for the pullback of $\beta.$

But then, dividing the integral on the right side of (\ref{SuffIt}) by that on the left, we obtain
\[
\int_{0}^{1}a^{v-1}(1-a)^{u-1}da = \beta(v,u) = \frac{\Gamma(v)\Gamma(u)}{\Gamma(v+u)}
\]
from which (\ref{SuffIt}) follows.

Now considering (\ref{It}), and supposing both relevant integrals to converge, we then find that
\begin{eqnarray*}
\int_{\gamma}\alpha \beta^{v+w-1}
&=&
\frac{(-1)^{v+w}}{\Gamma(v+w)}
\int_{1}^{0}
\left( \int_{1}^{t}\gamma^{*}\beta \right)^{v+w-1}\gamma^{*}\alpha \;\;\;\;\; \mbox{from Definition \ref{Cii}}
\\
&=&
\int_{1}^{0}\frac{-1}{\Gamma(v)\Gamma(w)}
\left[\int_{t}^{1}
\left(\int_{t}^{z}\gamma^{*}\beta\right)^{v-1}\left(\int_{z}^{1}\gamma^{*}\beta\right)^{w-1}\gamma^{*}\beta(z)\right]
\gamma^{*}\alpha(t)\;\;\;\mbox{by (\ref{SuffIt})}
\\
&=&
\int_{1}^{0} \frac{(-1)^{v+1}}{\Gamma(w)}
\left[ \int_{[1,t]} \left( \int_{z}^{1} \gamma^{*}\beta \right)^{w-1} \gamma^{*}\beta (\gamma^{*}\beta)^{v-1} 
\right] \gamma^{*}\alpha \;\; \mbox{from Definition \ref{Cii}}
\\
&=&
\frac{(-1)^{v+w}}{\Gamma(w)} \int_{[1,0]} \left( \int_{1}^{z}\gamma^{*}\beta\right)^{w-1}(\gamma^{*}\beta)^{v}\gamma^{*}\alpha \;\;\;\;\;\mbox{from (\ref{newdef})}
\\
&=&
\int_{[0,1]}
[\gamma^{*}\alpha (\gamma^{*}\beta)^{v}]
(\gamma^{*}\beta)^{w-1} \;\;\;\;\;\;\;\;\mbox{from (\ref{newdefn2})}
\\
&=&
\int_{\gamma}(\alpha \beta^{v})\beta^{w-1},
\end{eqnarray*}
so that (\ref{It}) holds as claimed. $\hfill \Box$

It would be interesting to give some geometric interpretation of these definitions along the lines of the above discussion involving simplices.

\subsection{The comultiplication formula}
\label{CoMult}

It will be convenient to introduce the following notation: 
\[
\int_{\gamma \rightarrow z}\beta:=\int_{0}^{z}\gamma^{*}\beta.
\] Then we have the 

{\thm{
\label{cOmUlt}
Suppose that $\alpha$ and $\beta$ are 1-forms on some manifold $M$ and $\gamma$ and $\delta$ are paths on $M$ for which 
\[
\int_{\gamma}\alpha \beta^{s} \; \; \mbox{and}\;\;
\int_{\delta}\alpha \beta^{s}
\]
both converge. Suppose also that 
\[
\left|
\int_{\delta^{-1}}\beta
\right|
>
\left|
\int_{\gamma^{-1} \rightarrow z}\beta
\right|,
\]
and that for any $z$ and sufficiently large $N$, 
\[
\sum_{n=0}^{N}\left( \matrix{s \cr n} \right)
\left(\int_{\gamma^{-1} \rightarrow z} \beta \right)^{n} \left( \int_{\delta^{-1}}\beta \right)^{-n}
\]
is dominated by the limit as $N \rightarrow \infty.$ Then
\[
\int_{\gamma \delta}\alpha\beta^{s} 
=
\int_{\delta}\alpha \beta^{s}
+
\sum_{n=0}^{\infty}\int_{\gamma}\alpha \beta^{n}\int_{\delta}\beta^{s-n}
\]
where we interpret 
\begin{equation}
\label{Weird}
\int_{\delta}\beta^{s-n}
\end{equation}
as \[
\left( \matrix{s \cr n} \right)
\cdot \frac{\int_{\delta}\beta^{s}}{\int_{\delta}\beta^{n}}
\]
whenever (\ref{Weird}) does not converge.
}}

{\bf Proof:}
\begin{eqnarray*}
\int_{\gamma \delta}\alpha\beta^{s}
&=&
\frac{-1}{\Gamma(s+1)}\int_{\delta^{-1}\gamma^{-1}}
\left( -\int_{(\delta^{-1}\gamma^{-1}) \rightarrow z}\beta\right)^{s}\alpha(z)
\\
&=&
\frac{-1}{\Gamma(s+1)}\int_{\delta^{-1}}
\left( -\int_{\delta^{-1}\gamma^{-1} \rightarrow y}\beta\right) ^{s}\alpha(y)
+
\frac{-1}{\Gamma(s+1)}\int_{\gamma^{-1}}
\left( -\int_{\delta^{-1}\gamma^{-1}\rightarrow x}\beta\right)^{s}\alpha(x)
\\
&=&
\frac{-1}{\Gamma(s+1)}\int_{\delta^{-1}}
\left( -\int_{\delta^{-1}\rightarrow y}\beta\right)^{s}\alpha(y)
+
\frac{-1}{\Gamma(s+1)}\int_{\gamma^{-1}}
\left( -\int_{\delta^{-1}}\beta - \int_{\gamma^{-1} \rightarrow x}\beta\right) ^{s}\alpha(x)
\\
&=&
\int_{\delta}\alpha\beta^{s} 
+ 
\frac{-1}{\Gamma(s+1)}\int_{\gamma^{-1}}
\left( -\int_{\delta^{-1}}\beta\right)^{s}\left(
\sum_{n=0}^{\infty}\left( \matrix{s\cr n} \right)
\left(-\int_{\delta^{-1}}\beta\right)^{-n}\left(
-\int_{\gamma^{-1} \rightarrow x}\beta\right)^{n}\right) \alpha(x)
\\
&=&
\int_{\delta}\alpha\beta^{s}+
\frac{-1}{\Gamma(s+1)}
\sum_{n=0}^{\infty}\left(\left( -\int_{\delta^{-1}}\beta\right)^{s-n}\right)
\cdot
\left( \matrix{s \cr n}\right)
\int_{ \gamma^{-1}}
\left( - \int_{\gamma^{-1} \rightarrow x}\beta \right)^{n}\alpha(x)
\\
&=&
\int_{\delta}\alpha\beta^{s}
+
\frac{1}{\Gamma(s+1)}
\sum_{n=0}^{\infty}\Gamma(s-n+1)
\int_{\delta}\beta^{s-n}
\cdot
\left( \matrix{s \cr n}\right)
n!\int_{\gamma} \alpha \beta^{n}
\\
&=&
\int_{\delta}\alpha\beta^{s}
+
\sum_{n=0}^{\infty}
\int_{\gamma}\alpha\beta^{n} \cdot \int_{\delta}\beta^{s-n}
\end{eqnarray*}
using the binomial series; and taking $z$ to be a parameter on $\delta^{-1}\gamma^{-1}$, with $y$ and $x$ parameters on $\delta^{-1}$ and $\gamma^{-1}$ respectively. $\hfill \Box$

With a view towards later application, we mention the  special case occuring when $\int_{\delta^{-1}}\beta = 0.$ It is immediate from the above proof that then,
\begin{equation}
\label{delta0}
\int_{\gamma \delta}\alpha \beta^{s-1} = \int_{\delta}\alpha \beta^{s-1}+\int_{\gamma}\alpha \beta^{s-1}.
\end{equation}

Of course, along with the coproduct formula, one would like some kind of product formula so as to have a {\sc Hopf} algebra of complex iterated integrals. It is obvious what such a formula would have to look like, but at this point there is a difficulty in the interpretation of the meaning of certain integrals in this formula, so the resolution of this problem will have to await further work.

Notice that in the coproduct formula, one only shifts by integers. For this reason, a less general iterative property than that discussed above probably suffices (i.e. in (\ref{It}) we would only need to consider pairs $(s,w)$ where one or other of the entries is an integer).

An important consequence of the comultiplication is the 
{\cor{
If $[\gamma]$ is a homotopy class of paths on $M$ for which there exists a representative $\gamma$ having the property that
\[
\left|
\int_{\gamma} \beta
\right|
\;>\;
\left|
\int_{\gamma^{-1}\rightarrow z}\beta \right|
\]
then for any path $\tilde{\gamma}$ in this homotopy class,
\[
\int_{\gamma}\alpha \beta^{s}
=
\int_{\tilde{\gamma}}\alpha \beta^{s}.
\]
}}
\newline {\bf Proof:}
Since $\int_{\gamma}\beta = \int_{\tilde{\gamma}}\beta$ by {\sc Cauchy's} theorem, 
the conditions of Theorem \ref{cOmUlt} apply to $\gamma^{-1} \tilde{\gamma}.$ Hence we have
\begin{eqnarray}
\int_{\gamma^{-1}\tilde{\gamma}}
\alpha \beta^{s}&
=&
\int_{\tilde{\gamma}}\alpha \beta^{s}
+
\sum_{n=0}^{\infty}
\int_{\gamma^{-1}}\alpha \beta^{n}
\cdot
\int_{\tilde{\gamma}}\beta^{s-n}\\
&=&
\label{homotloop}
\int_{\tilde{\gamma}}\alpha \beta^{s}
+
\sum_{n=0}^{\infty}
\int_{\gamma^{-1}}\alpha \beta^{n}
\cdot
\int_{{\gamma}}\beta^{s-n}
\end{eqnarray}
where the homotopy invariance of $\int_{\gamma}\beta^{t}$ -  which follows from the definition of this object in terms of $\int_{\gamma}\beta$ itself - gives (\ref{homotloop}).
Now since $\alpha$ and $\beta$ are holomorphic 1-forms on $M$ and $\gamma$ and $\tilde{\gamma}$ are homotopic, these forms are smooth over the region enclosed by the loop $\gamma^{-1} \tilde{\gamma}.$ But then by {\sc Cauchy's} theorem along with Definition \ref{Cii} it is evident that 
\begin{equation}
\label{hloop0}
\int_{\gamma^{-1} \tilde{\gamma}}
\alpha \beta^{s} =0.
\end{equation}

The same is true along the contractible $\gamma^{-1} \gamma$ - i.e.
\[
\int_{\gamma^{-1} \gamma} \alpha \beta^{s} = 0.
\]
Moreover, Theorem \ref{cOmUlt} applies to this path too, and we find
\begin{equation}
\label{contractloop}
0=\int_{\gamma^{-1}\tilde{\gamma}}
\alpha \beta^{s}
=
\int_{{\gamma}}\alpha \beta^{s}
+
\sum_{n=0}^{\infty}
\int_{\gamma^{-1}}\alpha \beta^{n}
\cdot
\int_{{\gamma}}\beta^{s-n}.
\end{equation}

Putting this all together by subtracting (\ref{contractloop}) from (\ref{homotloop}) married to (\ref{hloop0}), we find
\[
0 = \int_{\tilde{\gamma}}\alpha \beta^{s}
-
\int_{\gamma} \alpha \beta^{s}.
\]
$\hfill \Box$

\section{Integrating on $\Pomt$}

{\sc Deligne's} epic work [5] establishes the fundamental group of $\Pomt$ as an interesting object of study, providing as it does a test-case for the motivic philosophy. In work of {\sc Wojtkoviak} and {\sc Drinfel'd} related to this study, values of the {\sc Riemann} zeta function made a surprise appearance. This phenomenon is now well-understood: The  polyzeta numbers (also called multiple zeta values in the literature), 
are periods relating two distinct rational structures on the mixed {\sc Hodge} structure which comprises the {\sc Hodge} realization of the motivic fundamental group of $\lvP^{1} \backslash \{0,1,\infty\}$ with tangential base-point $\overrightarrow{01}$.
A fundamental reason for this is that  each polyzeta number admits an expression as an   iterated integral in the sense of {\sc Chen} over the holomorphic 1-forms of $\lvP^{1}\backslash \{0,1,\infty\}.$ By studying the complex iterated integrals defined above in the context of $\Pomt$, we are able to realize the polyzeta {\em functions} - along with other generalizations of the {\sc Riemann} zeta function - as iterated integrals along the tangential path from $0$ to $1.$

Preliminary to this investigation, 
we place a restriction on the functions $f$ which are being integrated
to ensure convergence of the integral under consideration:
{\defn{Let $k \in \lvZ_{\geq 0}.$ A $k$-{\sc Bieberbach} function is a function $f(z)$ which is holomorphic on the unit disk $D(0,1):=\{z \in \lvC :|z|<1\}$ and has a {\sc Taylor} series expansion
\[
f(z) = \sum_{n \geq 0}a_{n}z^{n}
\]
which satisfies the following property: $k$ is minimal for which there exist positive $N_{k}$ and $C_{k}$ so that
\[
|a_{n}| \leq C_{k}n^{{k}}
\]
whenever $n \geq N_{k}$ (i.e. $a_{n} = O(n^{k})$). 
\newline We shall say that a function is at least $k$-{\sc Bieberbach} if it is $l$-{\sc Bieberbach} for some $l \leq k.$}}

{\bf Examples: 
\newline 1.} Schlict functions are 1-{\sc Bieberbach} ({\sc De Branges}).
\newline {\bf 2.} $F_{\lvQ}(z):=\frac{z}{1-z}$ is 0-{\sc Bieberbach}. Later we show that $F_{\lvQ}$ underlies the {\sc Riemann} zeta function.

We point out that to make sense of an integral of the form of 
\[
\int_{0}^{t}\frac{dz}{z} 
\] 
some regularization of the logarithm at zero is necessary. In particular, should $f(z)$ be defined in some neighborhood $U$ of zero from which the points along the negative real axis have been deleted, so that for $\varepsilon$ close to zero and for $b \in U,$
\[
\int_{\varepsilon}^{b}f(z)dz
=
b_{0} + b_{1}\log \varepsilon + b_{2}(\log \varepsilon )^{2} + \ldots,
\] 
it is common usage to set
\[
\int_{0}^{b}f(z)dz
:=
b_{0}.
\]
By restricting the kinds of functions against which we integrate, we avoid this regularization issue altogether. To be precise, 
the $k$-{\sc Bieberbach} functions $f(z)$ henceforth considered will always have the property that $f(0) = 0,$ so that our complex iterated integrals will be at worst improper, and we can prove the

{\lemma{
\label{Lemmaconv}
Suppose $f(z)$ is at least $k$-{\sc Bieberbach} and vanishes at $z=0$. Then the improper integral
\[
\int_{0}^{1} \frac{(-\log z)^{s-1}}{\Gamma(s)}f(z)\frac{dz}{z}
=
\int_{[0,1]}f(z)\frac{dz}{z}\left(\frac{dz}{z}\right)^{s-1}
\]
converges for 
$Re (s) > k+1$.
}} 
\newline We shall henceforth use the notation
\[
\int_{[0,1]}f(z)\left(\frac{dz}{z}\right)^{s}
\]
for this integral.
\newline {\bf Proof:}
For $v \neq 0$ and  any $c \in (0,1],$
\[
\int_{0}^{c}\frac{(\log c-\log z)^{s-1}}{\Gamma(s)}z^{v}\frac{dz}{z} 
=
\frac{c^{v}}{v^{s}}
\]
via the  substitution $z^{v}=c^{v}u$ and use of the definition of $\Gamma(s).$ (Notice that like the integral defining $\Gamma(s),$ this integral is improper at $z=0.$)

The {\sc Taylor} series expression for $f$ converges uniformly on compacta in $D(0,1).$ Hence the order of the summation of this series and  integration over subintervals $[a,b]$ with $0<a <b <1$  may be interchanged.
Consequently, if
\[
f(z) = \sum_{n \geq 1} a_{n}z^{n}
\]
on $D(0,1),$ then for any $c \in (0,1),$
\begin{eqnarray*}
\int_{0}^{c} 
\frac{(\log c-\log z)^{s-1}}{\Gamma(s)}
f(z)\frac{dz}{z}
&=&
\lim_{\varepsilon \rightarrow 0}
\int_{\varepsilon}^{c}
\frac{(\log c - \log z)^{s-1}}{\Gamma(s)}
f(z)\frac{dz}{z}\\
&=&
\lim_{\varepsilon \rightarrow 0}
\sum_{n \geq 1}\int_{\varepsilon}^{c}
\frac{(\log c-\log z)^{s-1}}{\Gamma(s)}
a_{n}z^{n}\frac{dz}{z}\\
&=&
\sum_{n \geq 1}
\frac{a_{n}c^{n}}{n^{s}},
\end{eqnarray*}
and this sum converges by comparison with the sum for $\zeta(s-k)$ by  the $k$-{\sc Bieberbach} condition on the $a_{n}$.

For $c=1,$ the above argument does not suffice to allow for the interchange of the integral and sum. Instead, as in [15] we must resort to use of the {\sc Lebesgue} dominated convergence theorem since in $\Pomt$ the path of integration is not compact. The computation is most readily performed by means of exponential coordinates: To this end, let $-\log z=x,$ so that as improper integrals,
\[
\int_{0}^{1}\frac{(-\log z)^{s-1}}{\Gamma(s)}f(z)\frac{dz}{z}
=
\int_{0}^{\infty}\frac{x^{s-1}}{\Gamma(s)}f(e^{-x})dx
=
\int_{0}^{\infty}\frac{x^{s-1}}{\Gamma(s)}\sum_{n\geq 1}a_{n}e^{-nx}dx.
\]
Here the partial sums 
\[
\sum_{n=1}^{m}a_{n}e^{-nx}x^{s-1}
\]
are dominated by 
\[
F(k,\sigma,x):=\sum_{n=1}^{\infty}n^{k}e^{-nx}x^{\sigma-1}
\]
where $\sigma = \Re (s),$ making use of the $k$-{\sc Bieberbach} condition.
Now integrability of $F(k,s,x)$ over $[0,\infty)$ (provided that $\Re (s) >k+1$) suffices to complete the proof.

The only pole of $x^{1-s}F(k,s,x)$
along $[0,\infty)$ is at $x=0,$ but 
\[
{x^{k+1}}{\sum_{n=1}^{\infty}n^{k}e^{-nx}}
\]
is bounded on the unit disk centered at zero (say by $M_{0} e^{-\frac{1}{2}}$), since the limit of this expression at $x=0$ is finite:
Indeed, replacing $k$ in what follows by the least integer which exceeds it should $k \not \in \lvZ,$
\begin{eqnarray*}
{x^{k+1}}{\sum_{n=1}^{\infty}n^{k}e^{-nx}}&=&
x \sum_{n=1}^{\infty}(nx)^{k}e^{-nx}\\
&=&
x \left.\sum_{n=1}^{\infty}
\left( - \frac{d}{dt}\right)^{k}
e^{-tnx}\right|_{t=1} 
\\
&=&
x \left.\left( - \frac{d}{dt}\right)^{k} \frac{e^{-tx}}{1-e^{-tx}}\right|_{t=1},
\end{eqnarray*}
and one can compute the finite value of the limit as $x$ approaches 0 using {\sc L'H\^opital's} rule once the limit is interchanged with the $k$-fold derivative:
\[
\left.
\left(-\frac{d}{dt}\right)^{k} \lim_{x \rightarrow 0}
\frac{x}{e^{tw}-1}
\right|_{t=1}
=
\left.
\left(-\frac{d}{dt} \right)^{k} \frac{1}{t} \right|_{t=1}
= \left\{
\matrix{1 & {\mbox{if }} k=0\cr
        -(k-1)! & \mbox{if } k>0.}
\right.
\]

But then with 
\[
{x^{k+1}}{\sum_{n=1}^{\infty}n^{k}e^{-nx}}
\]
approaching zero as $x$ grows along $[0,\infty),$ for some $M_{1}$ also
\[
\left|
{x^{k+1}}{\sum_{n=1}^{\infty}n^{k}e^{-nx}}
\right| \leq M_{1}e^{-\frac{x}{2}}
\]
for all $x \geq 1.$ 
Let $M = \max \{M_{0}, M_{1}\},$ so that
\begin{eqnarray*}
\left|
\int_{0}^{\infty}
{x^{k+1}}{\sum_{n=1}^{\infty}n^{k}e^{-nx}}
x^{s-k-1}
\frac{dx}{x}
\right|
&\leq &
\int_{0}^{\infty}
\left|
{x^{k+1}}{\sum_{n=1}^{\infty}n^{k}e^{-nx}}
\right| x^{\sigma -k-1}\frac{dx}{x} \;\;\;
{\mbox{where $\sigma = \Re s$}}
\\
&\leq&
M\int_{0}^{\infty}
e^{-\frac{x}{2}}x^{\sigma-k -1} \frac{dx}{x}
\\
&=&
M 2^{\sigma-k -1} \Gamma(\sigma-k -1)\;\;\; 
\end{eqnarray*}
by again making use of the definition of the gamma function. Notice here that should $(\sigma - k-1) >0,$ this last  expression would be finite and $F(k,s,x)$ integrable. $\hfill \Box$

For the straight line path $[t,1]$ where $t  \in (0,1),$ the iterativity property was known classically, since in this case, the
definition accords with the classical fractional integral (see [1]): If $x$ is a real variable with $x \in (0, \infty)$ and $f(t)$ is integrable on this interval, then one notices that by repeated integration by parts, 
$\frac{1}{(k-1)!}\int_{0}^{x}(x-t)^{k-1}f(t)dt$ may be considered to be a $k$-fold integral of $f.$
For $\Re s >0,$ we then define the $s$-fold integral of $f$ as
\[
I_{s}f(x):=\frac{1}{\Gamma(s)}\int_{0}^{x}(x-t)^{s-1}f(t)dt.
\]
By the change of variables $-\log z = x-t,$ this may be seen to coincide with
\[
\int_{[e^{-x}, 1]}f(x-\log z) \frac{dz}{z} \left( \frac{dz}{z} \right)^{s-1}.
\]
Now the (classical) additivity property $I_{s}I_{t} = I_{s+t}$ of the operator $I_{s}$ amounts to our iterativity property.

On $\Pomt$ we may also take $\beta = \frac{dt}{1-t}$, and in integrating over the tangential path $[0,1]$, it is possible to show that such iterated integrals are characterized by the fact that they interpolate those integrals where the iteration occurs an integer number of times, while satisfying a suitable iterative property. We proceed to prove this.

Recall that (\ref{it}) is valid for any differential 1-form $\beta$, and as discussed before it has a unique interpolation  once 
a choice has been made of a branch of the logarithm. Therefore, if an iterative property can be established in the case of iteration over $\frac{dt}{1-t}$ along the path $[0,u]$ for $0<u\leq 1,$ necessarily
\[
\int_{[0,u]}\left(\frac{dt}{1-t}\right)^{s-1}
=
\frac{1}{\Gamma(s)}\left(\int_{[0,u]}\frac{dt}{1-t}\right)^{s-1}.
\] 
But then for any $k$-{\sc Bieberbach} $f$ for which $f(0)=0$,
\begin{eqnarray*}
\int_{[0,1]}f(z)\left(\frac{dz}{z}\right)^{s-1}
&=&
\int_{0}^{1}
\frac{(-\log z)^{s-1}}{\Gamma(s)}f(z)\frac{dz}{z}
\\
&=&
-\int_{1}^{0}\frac{(-\log (1-t))^{s-1}}{\Gamma(s)}f(1-t)
\frac{dt}{1-t} \;\;\;\;\;\;\; (t=1-z)\\
&=&
\int_{0}^{1}\left(\int_{[0,u]}\left(\frac{dt}{1-t}\right)^{s-1}\right)f(1-u)\frac{du}{1-u}.
\end{eqnarray*}
We shall use the notation $\int_{[0,1]}\left(\frac{dt}{1-t}\right)^{s-1}f(1-t)\frac{dt}{1-t}$ for this last integral expression, which is justified by the iterative property to follow. 

Recall that the iterativity property is the statement that
for a fixed $r \in \lvC$ with $\Re(r)>k+1,$ then  for 
any $w \in \lvC$ with $\Re(r)> \Re(w)>k+1$, 
and for any $g(z)$ for which $g(1-z)$ is 
$k$-{\sc Bieberbach} and has $g(1)=0$, it follows that
\begin{equation}
\label{IPA}
\int_{[0,1]}\left( \int_{[0,v]}\left(\frac{dt}{1-t}\right)^{r-w-1}  \right)
\left(\frac{dt}{1-t}\right)^{w}g(t)\frac{dt}{1-t}
=
\int_{[0,1]}
\left(
\frac{dt}{1-t} 
\right)^{r-1}
g(t)\frac{dt}{1-t}.
\end{equation}

Framed in a different way, 
\[
\int_{0}^{1}\int_{0}^{u}
\frac{(-\log(1-t))^{r-w-1}}{\Gamma(r-w)}
\frac{(\log (1-t)-\log (1-u))^{w-1}}{\Gamma(w)}
\frac{dt}{1-t}g(u)\frac{du}{1-u}
=
\int_{0}^{1}\frac{(-\log (1-u))^{r-1}}{\Gamma(r)}
g(u)\frac{du}{1-u}.
\]

Again this statement is a consequence of the non-trivial classical fact that the beta integral has an expression in terms of values of the gamma function: Indeed,
for
\begin{equation}
\label{beta0}
\frac{(-\log(1-u))^{r-1}}{\Gamma(r)}
=
\int_{0}^{u}
\frac{(-\log(1-t))^{r-w-1}}{\Gamma(r-w)}
\frac{(\log (1-t)-\log (1-u))^{w-1}}{\Gamma(w)}
\frac{dt}{1-t}
\end{equation}
to hold,  
\[
\beta(w,r-w) = \frac{\Gamma(w)\Gamma(r-w)}{\Gamma(r)}
= 
\int_{0}^{u}
\left( \frac{\log(1-t)}{\log (1-u)}\right)^{r-w-1}
\frac{(\log (1-t) - \log (1-u))^{w-1}}{(-\log (1-u))^{w}} \frac{dt}{1-t}
\]
must be true, and it {\em is} since the substitution
\[
\frac{\log (1-t)}{\log (1-u)} = y
\]
can be made to show that the integral is the same as
\[
\int_{0}^{1}y^{r-w-1}(1-y)^{w-1}dy = \beta(w,r-w).
\]

{\thm{For $(\Pomt, \frac{dt}{1-t}, [0,1])$, Defintion \ref{Cii}  is the only interpolation possible for which the iterativity property (\ref{IPA}) holds.
}}

{\bf Proof:}
The proof of the iterativity property works in this case since  the definition implies that
\begin{equation}
\label{woei}
G(u,s,w):=
\int_{0}^{u}\frac{(-\log (1-t))^{s-1}}{\Gamma(s)}\left(\frac{dt}{1-t}\right)^{w}
 = 
\int_{0}^{u}
\frac{(-\log (1-t))^{s-1}}{\Gamma(s)}
\frac{(\log (1-t)-\log (1-u))^{w-1}}{\Gamma(w)}\frac{dt}{1-t},
\end{equation}
where we write $s=r-w$. Should some other such integral expression exist, say
\[
G(u,s,w)= \int_{0}^{u}
\frac{(-\log (1-t))^{s-1}}{\Gamma(s)}F_{w}(t)
\frac{dt}{1-t},
\]
then for integer $w=n >1,$ in fact 
\[
F_{n}(t) =  
\frac{(\log (1-t)-\log (1-u))^{n-1}}{(n-1)!}
\]
because the usual antipode property may be used to unravel the iterated integral $G(u,s,n).$ For $F_{w}(t)$ to be a function in $w$ interpolating the $F_{n}(t),$ necessarily 
$$F_{w}(t) = e^{2\pi i r}\frac{(\log (1-t)-\log (1-u))^{w-1}}{\Gamma(w)}$$ by the  considerations pertaining to complex powers discussed before, for some integer $r$. Of course, here $r=0$ since there is no such exponential factor in the known expression for $G(u,s,w)$ in (\ref{woei}). 
$\hfill\Box$

\subsection{Multiplicative iterativity}

The above development of the complex iterated integral takes as departure point the iterative property which is necessarily satisfied. When iterating $\frac{dz}{z}$ over the tangential path $[0,1]$ in $\Pomt$, however, a second iterativity property holds.

In order to show this, we firstly prove an important fact which follows as a simple consequence of another computationally useful property, namely the power invariance of the iterated integral:
{\prop{{\bf{({\sc Haar} Property})}
\label{HProp}
\newline
Suppose that $f(z)$ is $k$-{\sc Bieberbach}.
Let $\alpha$ denote a positive real number. Then
\[
\int_{[0,1]}f(z^{\alpha})\left(\alpha \frac{dz}{z}\right)^{s}
=
\int_{[0,1]}f(z) \left( \frac{dz}{z} \right)^{s}.
\]
}}
As the notation suggests, 
\[
\int_{[0,1]}g(z)\left(\alpha \frac{dz}{z}\right)^{s}
:=
\int_{0}^{1}\frac{(-\alpha\log z)^{s-1}}{\Gamma(s)}g(z)\alpha \frac{dz}{z}
=\alpha ^{s} \int_{0}^{1}g(z)\left(\frac{dz}{z}\right)^{s}.
\]
Motivating the term `power invariance' is the fact that
\[
d \log z^{\alpha} = \frac{dz^{\alpha}}{z^{\alpha}}=
\frac{\alpha z^{\alpha -1}dz}{z^{\alpha}}
=\alpha \frac{dz}{z},
\]
and in the coordinates on $\lvC$ this amounts to invariance under the multiplicative group, hence the reference to {\sc Haar} measures.
\newline {\bf Proof:} Again it suffices to show the statement for $f(z) = z^{k}$ where $k$ is a non-negative integer. In this case, a direct computation involving a substitution $z^{\alpha k} = v$ shows the left side to equal 
\[
\frac{1}{k^{s}}.
\]
Then observe that this is the value of the right side of the equation we are proving, via the substitution $z^{k} = u.$
$\hfill \Box$

For future reference, we state an important fact implicit both in the proof of Lemma {\ref{Lemmaconv} and in that of Proposition \ref{HProp}, namely
{\cor{
\label{CorHaar}
For integer $k>0$ and $\Re s>1$, 
\[
\int_{[0,1]}z^{k}\left(\frac{dz}{z}\right)^{s} = \frac{1}{k^{s}}.
\]}}
\newline {\bf Proof:}
Take $f(z) = z$ and $\alpha = k$ in the {\sc Haar} property. Then use the fundamental fact, which is essentially equivalent to the definition of the Gamma function, that for any $s,$
\[
\int_{[0,1]}z\left(\frac{dz}{z}\right)^{s} = 1. 
\]
$\hfill \Box$

Now for notational ease make  the 
{\defn{If $F(z) = \sum_{n=-m}^{\infty}a_{n}z^{n}$ and $\Re s >1,$   define the $s$-gap transform of $F$ to be 
\[
^{s}F(z):=\sum_{n=-m}^{\infty}a_{n}z^{n^{s}}.
\]
}}

The following result furnishes an alternative means of defining complex iterated integrals of certain functions: 
{\thm{
If $F(z) = \sum_{n=1}^{\infty} a_{n}z^{n}$ is at least $k$-{\sc Bieberbach} for some $k\geq 0$ and $\Re s > k+1,$ then 
\[
\int_{[0,1]}F(z) \left( \frac{dz}{z} \right)^{s}
=
\int_{0}^{1} \; ^{s}F(z) \frac{dz}{z}.
\]
}} 
\newline {\bf Proof:}
From Corollary \ref{CorHaar}, when $n \neq 0$, 
\[
\int_{0}^{1}z^{n^{s}} \frac{dz}{z} = \frac{1}{n^{s}}.
\] 
But then invoking the interchange of the sum and integral by means of ideas as in the proof of Lemma \ref{Lemmaconv}, it follows that both sides of the equation give $\sum_{n=1}^{\infty}\frac{a_{n}}{n^{s}}. \hfill \Box$

Then we have:
{\thm{{\bf [Multiplicative Iterative Property]}
\label{wtshift}
If $F(z) = \sum_{n=1}^{\infty}a_{n}z^{n}$ is $r$-{\sc Bieberbach} for some $r>0,$ then for integer $k \geq 1,$ $^{k}F(z)$ is $\frac{r}{k}$-{\sc Bieberbach} and for $\Re s>r+k$,
\[
\int_{[0,1]}F(z) \left( \frac{dz}{z} \right)^{s}
=
\int_{[0,1]} \; ^{k}F(z) \left( \frac{dz}{z} \right)^{{s}/{k}}.
\]
}}
\newline {\bf Proof:}
That $^{k}F(z)$ is $\frac{r}{k}$-{\sc Bieberbach} is a triviality, and the equality of integrals follows from
\[
\int_{[0,1]}z^{n^{k}}\left(\frac{dz}{z} \right)^{{s}/{k}} = 
\frac{1}{(n^{k})^{{s}/{k}}} = \frac{1}{n^{s}},
\]
where we are again using Corollary \ref{CorHaar}. $\hfill \Box$

We remark that the failure of the expression for $^{s}F(z)$ given in the definition of the $s$-gap transform to be a power series for general non-integer $s$ precludes the proof of a more general iterative property without use of a vastly more complicated approach.

\subsection{The {\sc Riemann} zeta function as an iterated integral}

Using the formalism of the (additive) iterated integrals, {\sc Abel's} integral 
\begin{equation}
\label{Abel}
\zeta(s) = \frac{1}{\Gamma(s)}\int_{0}^{\infty}\frac{x^{s}}{e^{x}-1}\frac{dx}{x}
\end{equation}
may be thought of as an $s$-iterated integral:
{\thm{Whenever $\Re (s)>1,$ then
\label{Rzetait}
\begin{equation}
\zeta(s) = \int_{[0,1]}\frac{dz}{1-z}\left(\frac{dz}{z}\right)^{s-1}\; .
\end{equation}
}}
\newline 
{\bf Proof:} In (\ref{Abel}), make the substitution $-\log z = x$ and invoke the definition of the iterated integrals. 

Alternately, the statement is immediate from a direct computation using Corollary \ref{CorHaar} (with $k$ assuming successive positive integer values over which the sum is then taken).
$\hfill \Box$

A related family of complex iterated integrals for the {\sc Riemann} zeta function emerges from use of multiplicative iterativity:
\[
\int_{[0,1]}\sum_{n=1}^{\infty}z^{n^{k}} 
\left(\frac{dz}{z} \right) ^{{s}/{k}} 
=
\int_{[0,1]} 
\sum_{n=1}^{\infty}z^{n}\left( \frac{dz}{z} \right) ^{s} = \int_{[0,1]}\frac{z}{1-z}\left(\frac{dz}{z}\right)^{s}=
\zeta(s)
\]
for any integer $k \geq 1.$ Notice that when $k=2,$ we recover the integral for $\zeta(s)$ involving the theta function 
\[
\zeta(s) = \frac{1}{\pi^{{-s}/{2}}\Gamma(\frac{s}{2})}\int_{0}^{\infty}x^{{s}/{2}}\sum_{n=1}^{\infty}e^{-\pi xn^{2}}\frac{dx}{x},
\]
by the change of variables $z=e^{-\pi x}.$ The multiplicative iterativity of this integral expression contrasts with the additive iterativity of {\sc Riemann's} integral (\ref{Abel}). This is very interesting if one remembers that {\sc Riemann} gave two proofs of the functional equation for $\zeta(s)$, one using each of these integral expressions. However, he does not seem to have been aware of any complementarity in these perspectives, or if he was makes no mention thereof.

We remark that this discussion could be carried out with iteration over $\frac{dt}{1-t}$ instead. Then for example
\[
\zeta(s) = \int_{[0,1]}
\left( \frac{dt}{1-t} 
\right)^{s-1}\frac{dt}{t}.
\]

Including the local zeta function $\pi^{-\frac{s}{2}}\Gamma(\frac{s}{2}),$ (i.e. the (local) zeta function at the real prime), as a factor with $\zeta(s)$,
we obtain the function $Z (s)$ which may also be expressed as an $s$-iterated integral as follows:
{\thm{For $\Re s >1$,
\[
Z(s) = \int_{[0,1]}\frac{dx}{1-x}\left( \frac{1}{2\pi}\frac{(-\log x) 
dx}{x} \right)^{\frac{s-1}{2}}
\]}}

In the proof, we require the doubling formula of {\sc Legendre} for the factorial 
function:
\begin{equation}
\label{Legendredouble}
\Gamma(2z) = 2^{2z-1}\pi^{-{1}/{2}}\Gamma(z)\Gamma(z+\frac{1}{2}).
\end{equation}

{\bf Proof of Theorem:} 
\begin{eqnarray*}
Z(s) &=&
\pi^{-{s}/{2}}\Gamma(\frac{s}{2})\zeta(s)\\
&=&
\frac{\pi^{-{s}/{2}}\Gamma(\frac{s}{2})}{\Gamma(s)}\int_{0}^{1}\frac{(-\log 
x)^{s-1}dx}{1-x}\\
&=&
\frac{\pi^{{(1-s)}/{2}}2^{1-s}}{\Gamma(\frac{s-1}{2}+1)}
\int_{0}^{1}(-1)^{{(s-1)}/{2}}
\frac{(-(\log x)^{2})^{{(s-1)}/{2}}dx}{1-x}\;\;\;\;\;\;\;\;\;\;\;\;\;\mbox{by (\ref{Legendredouble})}\\
&=&
\frac{(-1)^{{(s-1)}/{2}}}
{2^{{(s-1)}/{2}}(2\pi)^{{(s-1)}/{2}}\Gamma(\frac{s-1}{2}+1)}
\int_{0}^{1} \left( \int_{1}^{x} -\frac{2\log x 
dx}{x}\right)^{{(s-1)}/{2}}\frac{dx}{1-x}\\
&=&
\int_{[0,1]}\frac{dx}{1-x}\left(-\frac{\log x dx}{2 \pi x} 
\right)^{{(s-1)}/{2}} 
\end{eqnarray*}
\hfill $\Box$

Recall that the functional equation of the {\sc Riemann} zeta function is the statement that this function $Z(s)$ is invariant under the transformation $s \mapsto 1-s.$

We remark that it is possible to develop {\sc Dirichlet} $L$-functions as iterated integrals using similar ideas. In particular, we find
that for a character $\chi$ of conductor $f$, 
\[
L(s, \chi)
= 
\int_{[0,1]} \frac{\sum_{a=1}^{f}\chi(a) z^{a}}{1-z^{f}} \left( \frac{dz}{z} \right)^{s},
\]
which may be developed using the multiplicativity property as
\[
L(s, \chi)
=
\int_{[0,1]}\sum_{a=1}^{f}\chi(a)\sum_{n=0}^{\infty}z^{(a+nf)^{k}}\left(\frac{dz}{z}\right)^{{s}/{k}}
\]
for any positive integer $k.$

\subsection{Multiple iterated integrals, polyzeta functions and polylogarithms}

The formalism may  be  extended to multiple versions of the iterated integrals by a simple induction argument based on the above. We perform this generalization in the case of the iteration of $\frac{dz}{z}$. 
For $j=1,2$ consider $k_{j}$-{\sc Bieberbach} functions $f_{j}(z)$ which are holomorphic at $z=0$ 
with also  $f_{j}(0)=0$.   Let $s_{j} \in \lvC$ have $\Re s_{j}>(k_{j}+j)$. Then
\[
\int_{[0,1]}f_{1}(z)\frac{dz}{z} \left(\frac{dz}{z}\right)^{s_{1}-1}
f_{2}(z)\frac{dz}{z} \left(\frac{dz}{z}\right)^{s_{2}-1},
\]
which is interpreted as
\[
\int_{[0,1]}\left[
\int_{[0,u]}f_{1}(z)\frac{dz}{z} \left(\frac{dx}{x}\right)^{s_{1}-1}
\right]
f_{2}(u)\frac{du}{u} \left(\frac{du}{u}\right)^{s_{2}-1}
\] 
converges by a similar argument to the one given before.{[{The vanishing of $f_{1}(z)$ at zero (so that the {\sc Taylor} series has first non-zero coefficient that of the linear term) facilitates the proof
since we can use the bound
\[
\left|
a_{1}b_{m-1}+\frac{a_{2}b_{m-2}}{2^{s_{1}}}
+ \ldots + \frac{a_{m}b_{0}}{m^{s_{1}}}\right|
\leq m C_{1}C_{2}m^{k_{2}}
\]
where the $a_{j}$ are coefficients for the power series for $f_{1}$ and the $b_{j}$ for the {\sc Taylor} series for  $f_{2}$; the $C_{1}$ factor bounds the coefficients $\left|\frac{a_{j}}{j^{s_{1}}}\right|$ with $0<j\leq m$ (taking $$C_{1} = \max\left\{\left|\frac{a_{j}}{j^{s_{1}}}\right|: 1 \leq j <N_{k_{1}}\right\} \cup \{C_{k_{1}}\}$$ with notation as in the definition of $k$-{\sc Bieberbach} functions),
and the $C_{2}m^{k_{2}}$ bound the $|b_{j}|.$ As seen before, the integral then converges provided 
$\Re s_{2} > (k_{2}+2).$}]}
Here, 
\[
\int_{[0,u]}
f_{1}(z)\frac{dz}{z} \left(\frac{dz}{z}\right)^{s_{1}-1}
=
\int_{0}^{u}
\frac{(\log u - \log z)^{s_{1}-1}}
{\Gamma{(s_{1})}}f_{1}(z)\frac{dz}{z}
=: h_{1}(u, s_{1})
\]
is a complex iterated  integral which satisfies an iterative property:
With notation as above (and polynomial $xp(x)$),
\begin{eqnarray*}
\int_{[0,u]}xp(x)\left(
\frac{dx}{x}\right)^{r}
&=&
\int_{0}^{u}\frac{(\log u - \log x)^{r-1}}
{\Gamma(r)}xp(x)\frac{dx}{x} \\
&=&
\int_{0}^{u}\int_{0}^{u_{1}}
\frac{(\log u-\log u_{1})^{r-w-1}}{\Gamma(r-w)}
\frac{(\log u_{1}-\log x)^{w-1}}{\Gamma(w)}
xp(x)\frac{dx}{x}\frac{du_{1}}{u_{1}}
\\
&=&
\int_{[0,u]}xp(x)\left(\frac{dx}{x}\right)^{w}
\left(\frac{du_{1}}{u_{1}}\right)^{r-w}
\end{eqnarray*}
follows from linearity by means of the substitution $v=\frac{x}{u}$ in the second expression, use of the iterativity property (\ref{It}), and then the substitution $u_{1} = u\tilde{u}$ for some intermediate variable $\tilde{u},$ followed by the reverse substitution $x=uv$. 

But then
\[
\int_{[0,1]}f_{1}(z)\frac{dz}{z} \left(\frac{dz}{z}\right)^{s_{1}-1}
f_{2}(z)\frac{dz}{z} 
\left(\frac{dz}{z}\right)^{s_{2}-1}
=
\int_{[0,1]}h(z, s_{1})f_{2}(z)\frac{dz}{z}
\left(
\frac{dz}{z}\right)^{s_{2}-1},
\]
where $h(z, s_{1})f_{2}(z)$ is $(k_{2}+1)$-{\sc Bieberbach} so that the iterative property holds not only in $s_{1},$ but in $s_{2}$ as well.

These ideas motivate the
{\defn{
\label{DMp}
Suppose that ${\bf{f}}:=(f_{1}(z), \ldots , f_{l}(z))$ is a tuple  of functions each holomorphic at $z=0$ with $f_{j}(0)=0$, such that 
$f_{j}(z)$ is $k_{j}$-{\sc Bieberbach}. 
Then the ${\bf{s}}:=(s_{1}, \ldots ,s_{l})$-multiple iterated integral of ${\bf{f}}$ against $\frac{dz}{z}$ is
\begin{eqnarray*}&&
\int_{[0,1]}f_{1}(z)\frac{dz}{z}\left(\frac{dz}{z}\right)^{s_{1}-1
}\cdots
f_{l}(z)\frac{dz}{z}\left(\frac{dz}{z}\right)^{s_{l}-1}\\
&:=&\int_{0}^{1} \int_{0}^{t_{l}}\cdots \int_{0}^{t_{2}}
\frac{(\log t_{2}-\log t_{1})^{s_{1}-1}}{\Gamma(s_{1})} f_{1}(t_{1})\frac{dt_{1}}{t_{1}} \cdots \\
&&
\frac{(\log t_{l}-\log t_{l-1})^{s_{l-1}-1}}{\Gamma(s_{l-1})} f_{l-1}(t_{l-1})\frac{dt_{l-1}}{t_{l-1}}
\frac{(-\log t_{l})^{s_{l}-1}}
{\Gamma(s_{l})} 
f_{l}(t_{l})\frac{dt_{l}}{t_{l}}
\end{eqnarray*}
provided that  $\Re (s_{j})>(k_{j}+j)$ 
for $1 \leq j \leq l.$
}}

Continuing the argument preceding the definition inductively, we establish the 
{\thm{{\bf{(Multiple Iterative Property)}}
For a fixed ${\bf{r}}=(r_{1}, \ldots, r_{l}) \in \lvC^{l}$ with $\Re (r_{j})>k_{j}+j$ 
for $1 \leq j \leq l$ then  for any 
$(s_{1}, \ldots, s_{l}) \in \lvC^{l}$ with 
$\Re (r_{j})>\Re (s_{j})>k_{j}+j$ 
for $1 \leq j \leq l$, 
and for any tuple 
$(f_{1}(z), \ldots, f_{l}(z))$ of functions each holomorphic at $z=0$, with $f_{j}(z)$ 
vanishing at $z=0,$ and with $f_{j}(z)$ being $k_{j}$-{\sc Bieberbach},
it follows that writing $w_{j} = r_{j}-s_{j},$ we have
\[
\int_{[0,1]}
f_{1}(z)\left(\frac{dz}{z}\right)^{s_{1}}
\left(\frac{dz}{z}\right)^{w_{1}}
\cdots
f_{l}(z)\left(\frac{dz}{z}\right)^{s_{l}}
\left(\frac{dz}{z}\right)^{w_{l}}
=\int_{[0,1]}f_{1}(z)
\left(\frac{dz}{z}\right)^{r_{1}}\cdots f_{l}(z)\left(\frac{dz}{z}\right)^{{r}_{l}}.
\]
}}
Otherwise stated, 
\[
\int_{0}^{1}\int_{0}^{\tilde{t}_{l}}\int_{0}^{t_{l}}
\cdots
\int_{0}^{t_{2}}\int_{0}^{\tilde{t}_{1}}
\frac{(\log \tilde{t}_{1}-\log t_{1})^{s_{1}-1}}{\Gamma(s_{1})}f_{1}(t_{1})
\frac{d{t_{1}}}{t_{1}}
\frac{(\log t_{2} - \log \tilde{t}_{1})^{w_{1}-1}}{\Gamma(w_{1})}
\frac{d\tilde{t}_{1}}{\tilde{t}_{1}}
\]
\[\cdots
\frac{(\log t_{l} - \log\tilde{t}_{l-1})^{w_{l-1}-1}}{\Gamma(w_{l-1})}
\frac{d\tilde{t}_{l-1}}{\tilde{t}_{l-1}}
\frac{(\log \tilde{t}_{l}-\log t_{l})^{s_{l}-1}}{\Gamma(s_{l})}
f_{l}(t_{l})\frac{d{t_{l}}}{t_{l}}
\frac{(-\log \tilde{t}_{l})^{w_{l}-1}}{\Gamma(w_{l})}\frac{d\tilde{t}_{l}}{\tilde{t}_{l}}
\]
\[
=
\int_{0}^{1}\int_{0}^{u_{l}}\cdots \int_{0}^{u_{2}}
\frac{(\log u_{2} - \log u_{1})^{r_{1}-1}}{\Gamma(r_{1})}
f_{1}(u_{1})\frac{du_{1}}{u_{1}}
\]
\[\cdots
\frac{(\log u_{l}-\log u_{l-1})^{r_{l-1}-1}}{\Gamma(r_{l-1})}
f_{l-1}(u_{l-1})\frac{d u_{l-1}}{u_{l-1}}
\frac{(-\log u_{l})^{r_{l}-1}}{\Gamma(r_{l})}
f_{l}(u_{l})\frac{du_{l}}{u_{l}}.
\]

Next we explain how the polyzeta functions may be expressed as complex iterated integrals.
For  integers $n_{j}$, it is well-known that the polyzeta numbers (also referred to as multiple zeta values in the literature), may be expressed as $(n_{1}+\ldots + n_{k})$-fold iterated integrals 
\[
\zeta(n_{1}, \ldots,n_{k})
=
\int_{[0,1]}\frac{dz}{1-z}\left(\frac{dz}{z}\right)^{n_{1}}
\frac{dz}{1-z}\left(\frac{dz}{z}\right)^{n_{2}}\ldots
\frac{dz}{1-z}\left(\frac{dz}{z}\right)^{n_{k}}.
\]
Once again, this expression also makes sense when the $n_{j}$ are replaced by non-integral complex numbers $s_{j}.$

{\sc Goncharov} and {\sc Kontsevich} found the following integral representation for the polyzeta functions:
\begin{equation}
\label{RI}
\zeta(s_{1}, s_{2}, \ldots, s_{l}) 
=
\frac{1}{\Gamma(s_{1})}\frac{1}{\Gamma(s_{2})}\ldots \frac{1}{\Gamma(s_{l})}
\int_{0}^{\infty}\cdots \int_{0}^{\infty} \int_{0}^{\infty}
\frac{t_{1}^{s_{1}-1}dt_{1}}{e^{t_{1}+t_{2}+\ldots +t_{l}}-1}
\frac{t_{2}^{s_{2}-1}dt_{2}}{e^{t_{2}+t_{3}+\ldots +t_{l}}-1}
\cdots
\frac{t_{l}^{s_{l}-1}dt_{l}}{e^{t_{l}-1}-1},
\end{equation}
valid provided $\Re (s_{l-j+1}+\ldots+s_{l})>j$ for $1 \leq j \leq l.$ When $l=1,$ the integral is the same as
the expression for {\sc Riemann's} zeta function found by {\sc Abel}. (See [15].)

It happens that using
the $\Pomt$ coordinates $x_{1}, \ldots, x_{k}$ determined via $x_{1}=e^{t_{1}+\ldots + t_{k}}$ and $x_{j+1}=e^{-t_{j}}x_{j}$ for $j=1, \ldots, k-1,$ the  integral is
\begin{equation}
\label{GKP1}
\int_{0}^{1}\int_{0}^{x_{k}}\cdots \int_{0}^{x_{2}}
\frac{(\log x_{2} - \log x_{1})^{s_{1}-1}}{\Gamma(s_{1})}
\frac{dx_{1}}{1-x_{1}}
\cdots
\frac{(\log x_{k} - \log x_{k-1})^{s_{k-1}-1}}{\Gamma(s_{k-1})}
\frac{dx_{k-1}}{1-x_{k-1}}
\frac{( - \log x_{k})^{s_{k}-1}}{\Gamma(s_{k})}
\frac{dx_{k}}{1-x_{k}}
\end{equation}
which may be regarded in an obvious way as a $k$-fold iterated integral along $[0,1]$ (in the sense of {\sc Chen})  generalizing {\sc Abel's} integral. But comparing this to Definition {\ref{DMp}} it is clear that in fact
\[
\zeta(s_{1},  \ldots, s_{l}) = \int_{[0, 1]}\frac{dz}{1-z}
\left( \frac{dz}{z} 
\right)^{s_{1}-1}\frac{dz}{1-z}\ldots \left(\frac{dz}{z}\right)^{s_{l}-1}
\]
whenever
$Re (s_{j})>j$ for all $1 \leq j \leq l.$

A striking duality exists: As before iteration over $\frac{dz}{1-z}$ could also be developed. This would give
\[
\zeta(s_{1},  \ldots, s_{l}) = \int_{[0, 1]}
\left( \frac{dz}{1-z} 
\right)^{s_{l}-1}\frac{dz}{z}\ldots \left(\frac{dz}{1-z}\right)^{s_{1}-1}\frac{dz}{z}.
\]

We remark that use of a similar change of coordinates for the integral expression known for the polylogarithm functions (see [3]) and use of these same ideas yields
\[
Li_{(s_{1}, \ldots, s_{l})}(t) = 
\int_{[0, t]}\frac{dz}{1-z}
\left( \frac{dz}{z} 
\right)^{s_{1}-1}\frac{dz}{1-z}\ldots \left(\frac{dz}{z}\right)^{s_{l}-1}
\]
for any $t \in [0,1]$, which also holds provided 
$Re (s_{j})>j$ for all $1 \leq j \leq l.$ The same formula is valid when $[0,t]$ is taken to indicate the straight line path from 0 to any $t$ in the open unit disc.

Also,  multiple versions of the {\sc Hurwitz} zeta functions may  be defined, and by similar considerations these satisfy
\begin{eqnarray*}
\zeta(s_{1}, \ldots, s_{l}; z)
&:=&
\sum_{0<n_{1}<\ldots<n_{l}}\frac{1}{(z+n_{1})^{s_{1}}
\ldots (z+n_{l})^{s_{l}}}
\\
&=&
\int_{[0, 1]}\frac{x^{z-1}dx}{1-x}
\left( \frac{dx}{x} 
\right)^{s_{1}-1}\frac{dx}{1-x}\left(\frac{dx}{x}\right)^{s_{2}-1}
\ldots \frac{dx}{1-x}\left(\frac{dx}{x}\right)^{s_{l}-1}
\end{eqnarray*}
whenever $Re (s_{j})>j$ for all $1 \leq j \leq l.$
Notice that $\zeta(s_{1}, \ldots, s_{l};1) = \zeta(s_{1}, \ldots, s_{l}).$

The  integral expressions for the polyzeta and {\sc Hurwitz} zeta functions  may  be thought of as   homotopy functionals evaluated along the homotopy class of the path $[0,1]$ in the fundamental groupoid consisting of homotopy classes of paths in $\Pomt$ from the tangential basepoint $\overrightarrow{01}$ to the tangential basepoint $\overrightarrow{10}$.
This particular path is very important, since it is identified with the {\sc Drinfel'd} associator $\Phi$ under the isomorphism of the unipotent completion of this fundamental groupoid with complex coefficients, with the group-like elements under comultiplication, of the completion of the free associative algebra generated by two symbols over $\lvC$ (i.e. the algebra of non-commuting power series in two variables, say $A_{0}$ and $A_{1}$, with complex coefficients).

\section{Further applications}

\subsection{Irrationality of {\sc Dedekind} zeta functions}
\label{sec4}

{\defn{When $F(z)$ is some $k$-{\sc Bieberbach} function vanishing at $z=0$, we shall call
\[
L(F)(s):=\int_{[0,1]}F(z)\left(\frac{dz}{z} \right)^{s}
\]
the $L$-function of $F.$
}}

Then the $L$-function of $F_{\lvQ}(z) = \frac{z}{1-z}$ is the {\sc Riemann} zeta function, and we notice immediately that  $L(F_{\lvQ})(s) = \zeta(s)$ has a simple pole at $s=1$ while 
$F_{\lvQ}(z)$ has a simple pole at $z=1.$

On the other hand, if $\chi$ is a non-trivial {\sc Dirichlet} character of conductor $f$, the $L$-function of 
\[
F_{\chi}(z):= 
\sum_{a=1}^{f}
\frac{\chi(a)z^{a}}{1-z^{f}}
\]
is  the {\sc Dirichlet} 
$L$-function 
$L(s, \chi) = \zeta(s,0;\chi)$
of  $\chi$, but in this case, $F_{\chi}(z)$ has no pole at $z=1$ and  $L(F_{\chi})(s)$ is non-singular at $s=1.$

Now $z$ is a coordinate on $\Pomt$ while $s$ describes $\lvC$, so there is {\em a priori} no connexion between them. For this reason, the correspondence between a pole of a function of $z$ and a pole of an associated function of $s$ in the case of {\sc Dirichlet} $L$-functions may appear  somewhat surprising. It turns out to be a consequence of the existence of the analytic continuation for the $L$-functions in the style of {\sc Riemann's} integral expression giving the analytic continuation of $\zeta(s)$. As such, this correspondence
holds quite generally:

{\thm{
\label{sz1}
Suppose that $F(z)$ 
is $k$-{\sc Bieberbach} for some $k$, vanishes at $a=0$, and is meromorphic in some neighborhood of $z=1$.
Then 
$L(F)(s)$ has a pole at $s=1$ if and only if $\frac{F(z)}{z}$ has a pole of non-zero residue at $z=1.$
Moreover, any pole of $L(F)(s)$ at $s=1$ is simple. When the pole of $F(z)$ at $z=1$ is also simple, the residue  agrees with that of $L(F)(s)$ at $s=1$. More generally, if $F(z) = \sum_{n \geq -m} a_{n}(z-1)^{n},$ the residue of 
$L(F)(s)$ at $s=1$ is $\sum_{n=-m}^{-1}(-1)^{1-n}a_{n}.$
}}
\newline
{\bf Proof:} 
There exists
$m \geq 1$ for which
$G_{m}(x):=x^{m}F(e^{-x})$ is regular at $x=0.$ 
We now fix $m$ as follows: If $F(z)$ is regular at $z=1,$ take $m=1.$ Otherwise, let $m \geq 1$ be minimal such that $G_{m}(0) \neq 0$ but $G_{m}(x)$ is regular at $x=0.$ 

Then define 
\[
H(F)(s) = \int_{C}{(-x)^{s-1}}F(e^{-x}) dx
\]
where $C$ once again denotes the {\sc Riemann} contour, and consider $H(F)$ for $\Re (s) >m+k+1.$
Denoting the part of $C$ which is a loop about zero by $\gamma_{0},$
\[
\int_{\gamma_{0}}
x^{m}F(e^{-x}) \frac{x^{s}}{x^{m}}\frac{dx}{x}
=0
\]
by {\sc Cauchy's} integral theorem. Consequently, as in the proof of Theorem {\ref{thm3}}, we find that
\[
H(F)(s) = (e^{i \pi s} - e^{- i \pi s})\Gamma(s) L(F)(s)
\]
for all  $s$ with $\Re (s) >m+k+1$
and hence on all of $\lvC.$

Then again as in the proof of Theorem {\ref{thm3}}, we have 
\[
L(F)(s)
=\frac{1}{\Gamma(s)2 i \sin (\pi s)}H(F)(s) = \frac{\Gamma(1-s)}{2 \pi i}H(F)(s).
\]

But
\[
H(F)(1) = \int_{\gamma_{0}}
F(e^{-x}){dx}
\]
because the integrals along the real axis cancel each other out. Now by the residue theorem the integral is non-zero exactly when $F(e^{-x})$ has a pole of non-zero residue at $x=0,$ which is precisely when $\frac{F(z)}{z}$ has a pole of non-zero residue at $z=1.$ Such are the instances in which
$L(F)(s)$ has a simple pole at $s=1$.

Since the residue of $\Gamma(1-s)$ at $s=1$ is 1, we also see that the residue of $L(F)$ at $s=1$ is 
\[
\frac{1}{2 \pi i} \cdot 2 \pi i 
\Res_{x=0}F(e^{-x})
=\frac{1}{2 \pi i}\int_{\gamma_{1}}
F(z)\frac{dz}{z}
\]
where $\gamma_{1}$ is a positively oriented loop about $z=1$. Using the power series expansion for $\frac{1}{z}$ at 
$z=1,$ the statement about the residues follows.
\hfill $\Box$

Suppose now that $K$ is a number field of degree $N$ over $\lvQ$ and $I$ denotes the set of non-zero integral ideals of $K$. Consider the {\sc Dedekind} zeta function $\zeta_{K}(s),$ which is known to have a simple pole at $s=1.$ 
This function also has a complex iterated integral expression, as may be seen from the 
{\lemma{The power series 
\[
F_{K}(z):=
\sum_{\dba \in I}z^{N(\dba)}
\] 
is at least $1$-{\sc Bieberbach}.
}}
\newline
{\bf Proof:}
Notice that if $\nu (n)$ denotes the number of ideals of $I$ of norm equal to $n$, we have
\[
\sum_{\dba \in I}z^{N(\dba)} = \sum_{n=1}^{\infty}\nu(n)z^{n}.
\]
Now $\sum_{j=1}^{m}\nu(j)  = \rho_{K} m + O(m^{1-\frac{1}{N}})$ where $\rho_{K}$ is the residue of 
$\zeta_{K}(s)$ at $s=1.$ (See [17], for example.) The
rough estimate
$\nu (n) \leq Cn$ then suffices to prove the lemma.
\hfill $\Box$

Convergence of the power series is uniform on compacta in the disc. Hence we may write
\[
\int_{[0,1]}
\sum_{\dba \in I}z^{N(\dba)}
\left( \frac{dz}{z} \right)^{s}
= 
\sum_{\dba \in I} \frac{1}{N(\dba)^{s}}
=:\zeta_{K}(s)
\]
for $\Re (s) >2.$

This suggests an archimedean analogue of the {\sc Iwasawa} algebra: The zeta function of a number field can be viewed as a power series in 
\[
\Lambda_{\infty}:=\lvZ[[T]].
\]
How far this analogy with {\sc Iwasawa} theory can be taken is an interesting question. One would like to see that properties of the $L$-function of a power series  are reflected in those of the power series itself.

For example,
the function 
\[
F_{pr}(z):=\sum_{p \;\;\mbox{prime}}z^{p}
\]
is not analytically continuable beyond the boundary of 
the unit disk, by the {\sc Fabry} gap theorem. This 
ought to reflect properties of 
\[
L(F_{pr})(s) = \sum_{p \;\;\mbox{prime}}\frac{1}{p^{s}}.
\] 

As another example, consider
the power series $\sum_{n =1}^{\infty} \mu(n)z^{n}$, which  is known to be singular at $z=1.$ (See [9].) The iterated integral of this function is $\frac{1}{\zeta(s)}$, so by  Theorem \ref{sz1} we see that $\sum_{n=1}^{\infty}\mu(n)z^{n}$ is not holomorphic on any punctured neighborhood of $z=1.$

Using the ideas of the next section, it will become clear that the rationality of the values of the {\sc Riemann} zeta function at negative integers is a direct consequence of the fact that 
$F_{\lvQ}(z) = \frac{z}{1-z}$
is rational. In line with the general philosophy that ``zeta functions should be rational'' we might expect that $F_{K}(z)$ would  also be rational, but this is not true in certain cases:

From Theorem \ref{sz1}, because $\zeta_{K}(s)$ has a pole at $s=1,$ we know that 
$F_{K}(z)$ is not regular at $z=1.$ Should $F_{K}(z)$ have a pole there, 
in the {\sc Laurent} series expansion for $F_{K}(z)$ at $z=1,$ a finite (alternating) sum of coefficients of $F_{K}(z)$  would be equal to  the residue of $\zeta_{K}(s)$ at $s=1$, which is known to be given by
\[
\rho_{K}=
\frac{2^{r_{1}}(2 \pi )^{r_{2}}R_{K}}{w\sqrt{|d_{K}|}}h_{K}
\]
where $r_{1}$ denotes the number of real embeddings of $K$; $2r_{2}$ the number of complex embeddings; $R_{K}$ is the regulator; $h_{K}$ the class number; $d_{K}$ the discriminant; and  $w$ the number of roots of unity in $K$.

Should $F_{K}(z)$ be a rational function, 
it would have an expression as a ratio of polynomials with integer coefficients by an elementary argument given in [2]: $F_{K}(z)$  is expressible as a power series with integer coefficents on the unit disk; so with notation $\nu(n)$ as above, if
\[
\sum_{n=0}^{m}p_{n}z^{n} \cdot \sum_{n=1}^{\infty}\nu(n)z^{n}
=
\sum_{n=0}^{l}q_{n}z^{n},
\]
then infinitely many linear equations with integer coefficients $\nu(n)$ arise, among which there is a solution in integers given that some solution exists. In this way, a  {\sc Laurent} series expansion about $z=1$ with rational coefficients would ensue. By Theorem \ref{sz1}, some linear combination of certain of these coefficients would have to equal $\rho_{K}$. (The only singularities of rational functions are poles.)
But
$\rho_{K}$ is expected to always be irrational, and whenever it is, $F_{K}(z)$ could not be rational.
(The difficulty in proving irrationality of $\rho_{K}$ lies in the fact that in general, both $\pi$ and logarithms of units appear in the formula for $\rho_{K}$ and it is not obvious that the product of these factors remains transcendental.)

By a theorem of {\sc Petersson} (see [2]), any power series with integer coefficients about zero having radius of convergence 1 is either  not analytically continuable beyond the boundary of the unit disk or gives rise to a rational function on $\lvC$; while a theorem of {\sc Fatou} in [9] asserts that such a power series yields a function which is either rational or non-algebraic.

Hence, using Theorem \ref{sz1} we find
{\thm{For a number field $K$ for which  $\rho_{K}$ is irrational, $F_{K}(z)$ is non-algebraic and noncontinuable outside of the unit disk.
}}

Consequently, we can also state the
{\cor{
\label{MainCor}
Irrationality of $\rho_{K}$ is an obstruction to the existence of a contour integral proof of the analytic continuation and functional equation for $\zeta_{K}(s)$ along the lines of {\sc Riemann's} first proof of the functional equation of $\zeta(s).$
}}
\newline This is evident from the fact that the contour of principal interest in such a proof would loop about $z=1.$

\subsection{Iterated integrals and derivatives}
{\label{sec2}}

{\sc Euler} conceived of an ingenious way to assign meaning
to the divergent infinite sum
\[
\sum_{n=1}^{\infty}n^{k}
\]
for $k \geq 1.$ The argument uses {\sc Abel} summation but ignores the divergence of the series being manipulated.{\footnote{Perhaps the most surprising fact in connexion with this argument is that it gives the same (correct) values of the {\sc Riemann} zeta function at negative integers, as a more rigorous approach does!}}

Now let $a \in \lvN$ have $a \geq 2$ and define 
$$\xi_{a}(n)= \left\{
\matrix{1 & \mbox{if}\;\; n \not \equiv 0\;\; (a) \cr     
1-a & \mbox{if}\;\; n  \equiv 0 \;\;(a). \cr}     \right.$$
Also let 
\[
\Psi(t)=\frac{\sum_{n=1}^{a}\xi_{a}(n)t^{n}}{1-t^{a}}.
\]
Using {\sc Euler's} ideas, {\sc Katz} produced the following generalization of his formula:
\begin{equation}
\label{EKcool}
\left( t \frac{d}{dt} \right)^{m} \Psi(t) |_{t=1} = 
(1-a^{m+1})\zeta(-m),\end{equation} for positive integers $m$.

Using the formalism of complex iterated integrals, it is not hard to see that also
\[
\int_{[0,1]}\Psi(t) \left(\frac{dt}{t} \right)^{s} = (1-a^{1-s})\zeta(s).
\] whenever $\Re s>1.$

This remarkable interplay between iterated derivatives and integrals holds quite generally:
{\thm{\label{thm3}
{\em{[{\sc Gel'fand - Shilov}]}}
If $F(t)
= \sum_{n=1}^{\infty} a_{n}t^{n}$ 
is holomorphic on the unit disk centered at $t=0,$ and is also analytic in some neighborhood of $t=1,$ then as a function of $s$,
\[
\int_{[0,1]} F(t) \left(\frac{dt}{t} \right)^{s}
\]
admits an analytic continuation which at 
negative integers $-k$ is given by
\[
\left( t\frac{d}{dt} \right) ^{k}F(t) |_{t=1}
\]
}}
\newline {\bf Proof:} Let $G(x):=F(e^{-x})$ and observe that this function is analytic in a neighborhood of $x=0.$
Consider 
\[
H(s):= \int_{C}(-x)^{s}G(x) \frac{dx}{x}
\]
where $C$ is the {\sc Riemann} contour from $+\infty$ to $0$ and back avoiding the positive real axis and looping around $0$ once in the positive direction.
Also define 
\begin{eqnarray*}
L(F)(s)&:=&\frac{1}{\Gamma(s)} \int_{0}^{\infty}x^{s}F(e^{-x})\frac{dx}{x}\\
&=&
\int_{[0,1]}F(t) \left(\frac{dt}{t}\right)^{s},
\end{eqnarray*}
which converges for $\Re (s) > k+1$ if $F$ is $k$-{\sc Bieberbach}.

Then we can show that $H(s) = 2i \sin (\pi s)\Gamma(s) L(F)(s)$: Indeed, suppose that $\Re s >k+1.$ Then on the first piece of the contour $C,$ (above the real axis) we know that $(-x)^{s} = e^{s \log x-i\pi s}$ whereas along the last piece of the contour (below the real axis)
$(-x)^{s} = e^{s\log x +i \pi s}$. Also, because $\Re s >1,$ the integrand is non-singular at zero, so as the radius of the loop about zero tends to zero, the value of the integral about this circular piece of $C$ also approaches zero. Then
\begin{eqnarray*}
H(s) &=& (-e^{-i\pi s}+e^{i \pi s}) \int_{0}^{\infty}x^{s}G(x) \frac{dx}{x}\\
&=&
2i \sin (\pi s) \Gamma (s) L(F)(s).
\end{eqnarray*}
The integral $H(s)$ converges for all complex $s,$ because $F(e^{-x})$ is a power series in $e^{-x}$ having no constant term, so that $F(e^{-x})$ dominates  $x^{s}$ as $x$ approaches infinity. Also,  the convergence is uniform on compacta so the function of $s$ determined by $H$ is complex analytic. Hence, using well-known identities satisfied by the $\Gamma$ function to write
\begin{equation}
\label{GenDir}
L(F)(s) = \frac{\Gamma(1-s)}{2 \pi  i} \int_{C}
(-x)^{s} G(x) \frac{dx}{x},
\end{equation}
we see that $L(F)(s)$ is a function defined and analytic at all points other than (possibly) the poles of
$\Gamma(1-s)$ - i.e. for $s \not \in \lvN \backslash \{0\}.$

From the convergence of $L(F)(s)$ on some right half-plane in $\lvC,$ 
we know then that the function has at most finitely many 
poles - to wit, at integers $0, 1, \leq k+1.$

Consequently it certainly makes sense to investigate the value 
of $L(F)(s)$ at negative integers, which we 
proceed to do:

$G(x) = F( e^{-x})$  is analytic in some neighborhood 
of $0 \in \lvC.$ Then write 
$G(x) = \sum_{m=0}^{\infty} b_{m}\frac{x^{m}}{m!}$.
On the pieces of the {\sc Riemann} contour lying above 
and below the real axis, we again have that 
$(-x)^{-k} = e^{-k\log x}e^{-i \pi k}$ and 
$(-x)^{k}=e^{-k\log x}e^{+i \pi k}$ respectively. Thus 
the integrals along these pieces are identical, although 
opposite in sign since the paths run in opposite 
directions. Hence
\begin{eqnarray*}
L(F)(-k) &=& 
\frac{\Gamma(1+k)}{2\pi i}
\int_{C} (-x)^{-k}G(x) \frac{dx}{x}\\
&=&\frac{\Gamma(1+k)}{2\pi i}\left(
\int_{+\infty}^{0}(-x)^{-k}G(x) \frac{dx}{x}
+
\int_{|x| = \delta}(-x)^{-k}G(x) \frac{dx}{x}
+\int_{0}^{+\infty}(-x)^{-k}G(x) \frac{dx}{x}\right)\\
&=&
\frac{k!}{2 \pi i}\int_{|x|=\delta}(-x)^{-k}
\sum_{m=0}^{\infty}b_{m}\frac{x^{m}}{m!} \frac{dx}{x}\\
&=&
\frac{(-1)^{k}k!}{2 \pi i}\sum_{m=0}^{\infty}\frac{b_{m}}{m!}
\int_{|x|=\delta}x^{m-k} \frac{dx}{x}\; 
\mbox{from  uniform convergence of the sum}\\
&=&
(-1)^{k}k!\sum_{m=0}^{\infty}\frac{b_{m}}{m!}
\frac{1}{2\pi}
\int_{0}^{2\pi} x^{m-k}d\theta\\
&=&
(-1)^{k}k!\sum_{m=0}^{\infty}\frac{b_{m}}{m!}
\frac{1}{2\pi}(2\pi \delta_{m,k})\\
&=&
(-1)^{k}b_{k}
\end{eqnarray*}

At the same time, 
\begin{eqnarray*}
\left( t \frac{d}{dt} \right)^{k} F(t)|_{t=1}&=&
\left(-\frac{d}{dx} \right)^{k}F(e^{-x})|_{x=0}\\
&=&
(-1)^{k} \left(\frac{d}{dx} \right)^{k}\sum_{m=0}^{\infty}b_{m}\frac{x^{m}}{m!}|_{x=0}\\
&=&
(-1)^{k}b_{k} = L(F)(-k)
\end{eqnarray*}
$\hfill \Box$

The above theorem was expressed by {\sc Gel'fand} and {\sc Shilov} in terms of
generalized functions - in particular,  
they show that the normalized distribution $\frac{x_{+}^{s-1}}{\Gamma(s)}$
satisfies
\[
\frac{x_{+}^{s-1}}{\Gamma(s)}|_{s=-n} = \delta^{n}(x)
\]
where $$\int_{0}^{\infty}\delta^{n}(x)\phi(x)dx = \phi^{n}(0)$$ for any test function $\phi.$
(See [10]I.\S 3.5). This is the same statement as that given above, under the co-ordinate 
change $x = -\log t.$

The proof of Theorem {\ref{thm3}} may easily be modified to show 
{\thm{
\label{3}
For $F$ as above and $w \in (0,1)$ arbitrary, then the function
\[
\int_{[w,1]} F(t) \left(\frac{dt}{t} \right)^{s}
\]
has the {\em same} analytic continuation to negative integers as does 
\[
\int_{[0,1]} F(t) \left(\frac{dt}{t} \right)^{s}.
\]}}
This $w$-independence is quite surprising.
From the distribution viewpoint, it is certainly true that the analytic continuation at negative integers is some kind of derivative {\sc Dirac} distribution centered at zero (corresponding to $1 \in \lvP^{1} \backslash \{0,1, \infty\}$), but for $Re (s) >1$ the distribution is not even compactly supported! In the homotopy theory perspective, the natural notion of tangential base-point is in evidence here: The analytic continuation of the iterated integrals is the same for all paths which lie along the tangential path between 0 and 1 in $\Pomt$, which end in the tangential base point $\stackrel{\rightarrow}{01};$ but there is no apparent reason why this should be so and as a function of $w$ the iterated integral is certainly non-constant. Observe that this implies that the $p$-adic $L$-functions interpolate values of a {\em family} of functions at negative integers.


Theorem \ref{thm3} may be used to immediately write down the formula for the values of $L$-functions at negative integers. In particular, 
if $\chi$ is a non-trivial {\sc Dirichlet} character of conductor $f$, then 
\[
L(-m, \chi) = \left( t\frac{ d}{dt}\right)^{m} \left. \frac{\sum_{a=1}^{f}\chi(a)t^{a}}{1-t^{f}}\right|_{t=1}.\]

Moreover, for such a character $\chi$, since $\sum_{a=1}^{f}\chi(a) = 0$ we may use Theorem {\ref{thm3}} to see that the generalized {\sc Hurwitz} zeta function
\[
\zeta(s,z; \chi) = \int_{[0,1]}t^{z}
\sum_{a=1}^{f}\frac{\chi(a)t^{a}}{1-t^{f}}
\left(\frac{dt}{t} \right)^{s}
\]
has analytic continuation to the negative integers given by
\[
\left.
\left( t \frac{d}{dt} \right)^{n}
\left( t^{z}
\sum_{a=1}^{f}\frac{\chi(a)t^{a}}{1-t^{f}}
\right)\right|_{t=1}
=
-\frac{B_{n+1, \chi}(z)}{n+1}.
\]
Effectively this is  a rewriting of the definition of the generalized {\sc Bernoulli} polynomials using the generating series under the change of coordinates $t = e^{-w}.$ 

\subsection{Monodromy of polylogarithms}
\label{Mop}

Denoting the straight line path from $0$ to $w \in D'(0,1)$ by $[0 \rightarrow w]$, then as asserted before,
\[
Li_{s}(w) = \int_{[0 \rightarrow w]}\frac{dx}{1-x}\left( \frac{dx}{x} 
\right)^{s-1}.\]

This is because
\begin{eqnarray*}
\int_{[0 \rightarrow z]}x^{k}dx \left( \frac{dx}{x} \right)^{s-1} 
&= &
\int_{[0,1]}(zt)^{k}zdt \left( \frac{zdt}{zt} 
\right)^{s-1}\\
&=&z^{k+1}\int_{[0, 1]} t^{k}dt\left(\frac{dt}{t} 
\right)^{s-1}\\
&=& \frac{z^{k+1}}{(k+1)^{s}}
\end{eqnarray*}
Or seen another way:
\begin{eqnarray*}
\int_{[0 \rightarrow z]}x^{k}dx \left( \frac{dx}{x}\right)^{s-1}
&=&\frac{(-1)^{s-1}}{\Gamma(s)}\int_{0}^{z}\left(\int_{0}^{w} 
\frac{dx}{x} \right)^{s-1}w^{k}dw\\
&=&\frac{(-1)^{s-1}}{\Gamma(s)}\int_{0}^{1}
\left(\int_{0}^{r} \frac{dt}{t} \right)^{s-1}z^{k+1}t^{k}dt 
\;\; \mbox{using}\;\; w=zt\\
&=&z^{k+1} \int_{[0, 1]}t^{k}dt \left( \frac{dt}{t} 
\right)^{s-1}\\
&=& \frac{z^{k+1}}{(k+1)^{s}}.
\end{eqnarray*}

The sum of such terms for $k=1, 2, \ldots$ and the interchange of the sum 
and integral 
(which as before is allowed because of an argument involving use of the {\sc Lebesgue} dominated convergence theorem) results in the above formula for the polylogarithm functions.

A well-known fact with an elegant expression in terms of iterated integrals is the general monodromy theorem:
{\thm{\label{mop}
$Li_{s}(w)$ continued analytically along a loop 
$\gamma$ about 1 (i.e. 
the monodromy of the general polylogarithm function) is 
\begin{equation}
\label{monodpolylog}
Li_{s}(w) - \frac{2 \pi i}{\Gamma(s)} \log^{s-1}(w) 
=
\int_{[0 \rightarrow w]} \frac{dx}{1-x} \left(\frac{dx}{x} \right)^{s-1} 
+ 
\int_{\gamma}\frac{dx}{1-x} \cdot \int_{[1 
\rightarrow w]} \left( \frac{dx}{x} \right) ^{s-1}.\end{equation}
}}

Classically,  {\sc Jonqui{\`e}re's} formula was used to effect the proof, but a direct topological proof using the homotopy functional (iterated integral) perspective is desirable since (\ref{monodpolylog}) is reminiscent of a coproduct formula in which many terms vanish. We proceed to give such a proof:

{\bf Proof:}
Fix $w$ bounded away from 1 and let $\eta>0$ have 
\[
|\log \eta|<|\log w|.
\]
Now notice that $[0 \rightarrow w]$ is homotopic to
the composition of the straight 
line paths  in $\lvP^{1} \backslash \{0,1,\infty\}$ from 0
to $\eta$ and from $\eta$ to $w$, which  will be denoted 
$[0 \rightarrow \eta]$ and 
$[\eta \rightarrow w]$ respectively.
Consequently, the homotopy functional 
\[
\int_{\cdot} \frac{dx}{1-x} \left( \frac{dx}{x} \right )^{s-1}
\]
evaluated along $[0 \rightarrow w]$ gives the same value as when it is evaluated 
along the succession of paths $[0 \rightarrow \eta] \cdot [\eta\rightarrow w]$.
i.e. 
\begin{equation}
\label{fr}
Li_{s}(w) 
= 
\int_{[0 \rightarrow w]} \frac{dx}{1-x} 
\left( 
\frac{dx}{x} 
\right )^{s-1}=
\int_{[0 \rightarrow \eta] \cdot [\eta\rightarrow w]} 
\frac{dx}{1-x} \left( \frac{dx}{x} \right )^{s-1}.
\end{equation}
By the definition of the iterated integral, using also the proof - and notation - of Theorem \ref{cOmUlt}, it then follows that
\begin{equation}
\label{frog}
Li_{s}(w) = \int_{[0\rightarrow \eta]}\frac{1}{\Gamma(s)}\left(-\int_{[w \rightarrow \eta]\cdot [\eta \rightarrow 0] \rightarrow x} \frac{dz}{z}\right)^{s-1}\frac{dx}{1-x}
+
\int_{[\eta \rightarrow w]}\frac{dx}{1-x} \left( \frac{dx}{x} \right 
)^{s-1}.
\end{equation}

Now let $\varepsilon >0$ be fixed, and 
denote the  circular path of radius $\varepsilon$ about 1 by 
$\gamma_{1, \varepsilon}. $
Again invoking the homotopy invariance,  
the
analytic continuation of 
$Li_{s}(w)$ around 1 may now be expressed by
\[
\int_{[0\rightarrow \eta]}\frac{1}{\Gamma(s)}\left(-\int_{[w \rightarrow \eta]\cdot [\eta \rightarrow 0] \rightarrow x} \frac{dz}{z}\right)^{s-1}\frac{dx}{1-x}
+
\int_{[\eta \rightarrow 1-\varepsilon] \cdot 
\gamma_{1, \varepsilon} \cdot 
[1-\varepsilon \rightarrow w]} \frac{dx}{1-x} \left( \frac{dx}{x} \right 
)^{s-1}, 
\]
the second term of which can be calculated using the coproduct formula of Theorem \ref{cOmUlt} applied associatively to the three paths along which the 
homotopy functional is being computed, since 
the choice of $\eta$ facilitates the following technical condition:
\[
\left|\int_{[w \rightarrow 1- \varepsilon]}\frac{dx}{x}\right|
\;
>
\;
\left| \int_{\gamma_{1, \varepsilon}^{-1} \cdot [1-\varepsilon \rightarrow \eta]\rightarrow z}\frac{dx}{x}\right|\]
for any $z \in [0,1]$ (again with notation as in Theorem \ref{cOmUlt}).

Now regarding the first two paths $[\eta \rightarrow 1-\varepsilon]$ and $\gamma_{1,\varepsilon},$ as per (\ref{delta0}), the coproduct in fact degenerates into the sum
\[
\int_{[\eta \rightarrow 1-\varepsilon]}\frac{dx}{1-x}\left(\frac{dx}{x}\right)^{s-1} +\int_{\gamma_{1,\varepsilon}}\frac{dx}{1-x}\left(\frac{dx}{x}\right)^{s-1}
\]
since $\int_{\gamma_{1, \varepsilon}}\frac{dx}{x}= 0.$

Hence, we obtain
\begin{eqnarray*}
& &
\int_{[\eta \rightarrow 1-\varepsilon] \cdot \gamma_{1,\varepsilon}\cdot [1-\varepsilon \rightarrow w]}
\frac{dx}{1-x} \left(\frac{dx}{x} \right)^{s-1}
\\
&=&
\int_{[1-\varepsilon \rightarrow w]}\frac{dx}{1-x} \left(\frac{dx}{x} \right)^{s-1}
+\sum_{n=0}^{\infty}\int_{[\eta \rightarrow 1-\varepsilon] \cdot \gamma_{1,\varepsilon}}\frac{dx}{1-x}\left(\frac{dx}{x}\right)^{n} \cdot \int_{[1-\varepsilon \rightarrow w]}\left(\frac{dx}{x}\right) ^{s-1-n}
\\
&=&
\int_{[1-\varepsilon \rightarrow w]}\frac{dx}{1-x} \left(\frac{dx}{x} \right)^{s-1}
+
\sum_{n=0}^{\infty}\int_{[\eta \rightarrow 1-\varepsilon]}
\frac{dx}{1-x}\left(\frac{dx}{x}\right)^{n} \cdot \int_{[1-\varepsilon \rightarrow w]}\left(\frac{dx}{x}\right) ^{s-1-n}
\\
& &
+
\sum_{n=0}^{\infty}\int_{\gamma_{1,\varepsilon}}\frac{dx}{1-x}\left(\frac{dx}{x}\right)^{n} \cdot \int_{[1-\varepsilon \rightarrow w]}\left(\frac{dx}{x}\right) ^{s-1-n}
\\
&=&
\int_{[\eta \rightarrow w]}\frac{dx}{1-x} \left(\frac{dx}{x} \right)^{s-1}
+
\sum_{n=0}^{\infty}\int_{\gamma_{1,\varepsilon}}\frac{dx}{1-x}\left(\frac{dx}{x}\right)^{n} \cdot \int_{[1-\varepsilon \rightarrow w]}\left(\frac{dx}{x}\right) ^{s-1-n},
\end{eqnarray*}
where the first term of the last expression results from  the comultiplication formula applied to the succession of paths $[\eta \rightarrow 1-\varepsilon] \cdot [1-\varepsilon \rightarrow w]$. Now in the remaining sum, allowing $\varepsilon$ to approach 0, the only non-zero integral about 
$\gamma_{1, \varepsilon}$ is the one for which $n=0$. This may be seen by writing the iterated integral as a contour integral, and recalling that here, $\Re s >1.$

Consequently, reintroducing the first term of (\ref{frog}), we end up with
\[
Li_{s}(w)+\lim_{\varepsilon \rightarrow 0}\int_{\gamma_{1, \varepsilon}}\frac{dx}{1-x} \cdot
\int_{[1-\varepsilon \rightarrow w]}\frac{dx}{x}^{s-1}
=
Li_{s}(w)-2\pi i \frac{(\log w)^{s-1}}{\Gamma(s)}.
\]
$\hfill \Box$

\parskip 3mm

{\bf{References}}
\parskip 0mm

[1] Andrews, G.; Askey, R.; Roy, R. {\em Special Functions}, Encyclopedia of Mathematics and its Applications, Cambridge University Press, Cambridge, 1999.

[2] Bieberbach, L. {\em Analytische Fortsetzung}, Ergebnisse der Mathematik, Springer Verlag, Berlin, 1955.

[3] Cartier, P. {\em Fonctions polylogarithmes, nombres polyz{\^e}tas et groupes pro-unipotents}, 
S{\'e}m. Bourbaki 2000-2001, exp. ${\mbox{n}}^{\circ}$
885.

[4] Chen, K.-T. {\em Iterated path integrals}, Bull. Amer. Math. Soc. {\bf 83} (1977), 831-879.

[5] Deligne, P. {\em Le groupe fondamental de la droite projective moins trois points}, in ``Galois groups over $\lvQ$: proceedings of a workshop held March 23-27, 1987'', editors K.Ribet et al. Springer-Verlag, New York, 1989, 79-297.

[6] Deligne, P.; Ribet, K. {\em Values of abelian $L$-functions at negative integers over totally real fields}, Inventiones math. {\bf{59}} (1980), 227-286.

[7] Drinfel'd, V.G. {\em On quasitriangular quasi-Hopf 
algebras and a group closely connected with $\Gal (\overline{\lvQ}/\lvQ)$}, Leningrad Math J. {\bf{2}} (1991) no.4, 829-860.

[8] Edwards, H.M. {\em Riemann's zeta function}, Dover, NY, 1974.

[9] Fatou, P. {\em S{\'e}ries trigonom{\'e}triques
et s{\'e}ries de Taylor}, Acta math. {\bf 30}, 1906, 335-400.

[10] Gel'fand, I.M.; Schilow, G.E. {\em Verallgemeinerte funktionen Teil I}: Verallgemeinerte Funktionen und das Rechnen mit ihnen.  Hochschulb\"ucher f\"ur Mathematik, Bd. 47 VEB Deutscher Verlag der Wissenschaften, Berlin 1960.

[11] Goncharov, A.B. {\em Multiple polylogarithms and mixed Tate motives}, arXiv: math.AG/0103059, 2001.

[12] Hain, R. {\em Classical Polylogarithms}, Motives (Seattle, WA, 1991), Proc. Sympos. Pure Math., {\bf 55}, Part 2, Amer. Math. Soc., Providence, RI, 1994, 3-42.

[13] Hain, R. {\em Periods of limit mixed Hodge structures}, Current developments in mathematics, Int. Press, Somerville, MA, 2003, 113-133.

[14] Hain, R. {\em Lectures on the Hodge-De Rham Theory of $\pi_{1}(\Pomt)$}, Arizona Winter School 2003.

[15] Hida, H. {\em  Elementary theory of $L$-functions and  Eisenstein series}, Cambridge University Press, Cambridge, 1993.

[16] Ihara, Y. {\em Braids, Galois groups, and some arithmetic functions}, Proceedings of the ICM, Kyoto, 1990, 99-120.

[17] Lang, S. {\em Algebraic Number Theory}, (Second Edition), Springer-Verlag, New York, 1994.

[18] Ngoc Minh, H.; Petitot, M.; van der Hoeven, J. {\em Shuffle algebra and polylogarithms}, Discrete Mathematics {\bf{225}} (2000), 217-230.

[19] Riemann, B. {\em Ueber die Anzahl der Primzahlen unter einer gegebenen Gr{\"o}sse}, in ``Gesammelte Werke'', Teubner, Leipzig 1892.

\end{document}